\documentclass[a4paper,11pt]{amsart}
\usepackage[latin1]{inputenc}

\usepackage{multicol}
\usepackage{amssymb, amsmath, amsthm}
\usepackage{latexsym}
\usepackage{dsfont}
\usepackage{graphicx}

\usepackage{cancel}
\usepackage[normalem]{ulem}

\usepackage{curves} 
\usepackage{xcolor}
 
\input xy
\xyoption{all}

\usepackage{verbatim}

\usepackage{hyperref}

\newtheorem{theorem}{Theorem}[section]

\newtheorem{corollary}[theorem]{Corollary}

\newtheorem{lemma}[theorem]{Lemma}

\newtheorem{proposition}[theorem]{Proposition}

	{\par\noindent{\bf Proposition \ref{res:hiper}.}\!\!
	\nopagebreak\normalsize\it}%

	{\par\noindent{\bf Theorem \ref{result43}.}\!\!
	\nopagebreak\normalsize\it}%

	{\par\noindent{\it Sketch of the proof}.  
	\nopagebreak\normalsize}%
	{\hfill\linebreak[2]\hspace*{\fill}$\circlearrowleft$}

	{\par\noindent{\it Proof of Proposition }\ref{prop:stab:smc}.  
	\nopagebreak\normalsize}%
	{\hfill\linebreak[2]\hspace*{\fill}$\circlearrowleft$}

	{\par\noindent{\it Proof of Propositions }\ref{adap:mon}{\it and }\ref{simult:adap}.\!\!\!
	\nopagebreak\normalsize}%
	{\hfill\linebreak[2]\hspace*{\fill}$\circlearrowleft$}

{\theoremstyle{definition}

       \newtheorem{remark}[theorem]{Remark}

       \newtheorem{parrafo}[theorem]{{\!}}  }

\numberwithin{equation}{theorem}

\newcommand{\calo}{{\mathcal {O}}}

\DeclareMathOperator{\Max}{\underline{Max}}
\DeclareMathOperator{\mult}{mult}

\DeclareMathOperator{\Proj}{Proj}

\DeclareMathOperator{\Spec}{Spec}

\newcommand{\M}{{\mathfrak M}}

\newcommand{\p}{{\mathfrak p}}

\newcommand{\m}{{\mathcal{M}}}
\renewcommand{\P}{\mathcal{P}}

\newcommand{\x}{{\mathbf{x}}}

\definecolor{darkpurple}{rgb}{0.28,0.24,0.55}
\definecolor{lightblue}{rgb}{0,0.75,1}

\title{On the simplification of singularities 
by blowing up at equimultiple centers}
\author{Orlando E. Villamayor U.}
\thanks{2000 {\em Mathematics subject classification. 14E15}}
\thanks{The author is  partially supported by MTM2009-07291}
\date{December 2015}
\address{Dpto. Matem\'aticas,  Universidad
Aut\'onoma de Madrid and Instituto de Ciencias Matem\'aticas CSIC-UAM-UC3M-UCM \\
Ciudad Universitaria de Cantoblanco, 28049 Madrid, Spain}
\email[Orlando E. Villamayor U.]{villamayor@uam.es}
\keywords{Multiplicity, integral closure of ideals, Rees rings.}
\begin{document}
\maketitle

\begin{abstract} Resolution of singularities of varieties over fields of characteristic zero can be proved by using the multiplicity as main invariant. The proof of this result leads to new questions in positive characteristic. We discuss here results which follow by induction on the dimension of the varieties. 

Fix a variety $X^{(d)}$ of dimension $d$ over a {\em perfect field} $k$ or, more generally, a pure dimensional scheme of finite type over $k$.
Fix a closed point $x\in X^{(d)}$ of multiplicity $e>1$. Define a {\em local simplification of the multiplicity at} $x\in X^{(d)}$ as a proper birational map, say
$X^{(d)}\leftarrow X^{(d)}_1$, where $X^{(d)}$ denotes now a neighborhood of $x$,  so that $X^{(d)}_1$ has multiplicity $<e$ at any point $x_1\in X^{(d)}_1$. 

Assume, by induction on $d$, the existence of local simplifications of the multiplicity for schemes over $k$ of dimension $d'$, for all $d' <d$. We prove, under this inductive assumption, that a local simplification at $x\in X^{(d)}$  can be constructed when $(C_{X,x})_{red}$ is not regular
Here $C_{X,x}$ denotes the tangent cone of $x\in X$, and $(C_{X,x})_{red}$ is the reduced scheme. The paper uses classical results of commutative algebra, and compares the effect of blowing up along equimultiple centers, and along normally flat centers.
\end{abstract}
{\tableofcontents}
\begin{section}{Introduction}
\begin{parrafo}

In this paper, which is largely expository, we discuss two invariants that are attached to a singular point of a variety. One is the multiplicity, and another is the Hilbert Samuel series at such point. 
 This leads to two stratifications on the variety into points with the same invariant. When the variety is a hypersurface both stratifications coincide, but in general they are different. So, in general, different properties hold if we blow up at centers 
where the multiplicity is constant, and when we blow up at centers where Hilbert Samuel invariant is constant. This leads to the notion of {\em equimultiple centers} and 
{\em normally flat centers} respectively, to be discussed along these notes (see also \cite{Herrmann}). Hironaka proves resolution of singularities in characteristic zero by using invariants that grow from the Hilbert Samuel function, and by blowing up at normally flat centers (\cite{Hironaka64}). A similar statement holds using the multiplicity, and blowing up at equimultiple 
centers (\cite{multiplicity}).

Fix a variety, or an excellent pure dimensional scheme $X$, and 
consider the multiplicity of the local ring $\calo_{X,x}$ at each $x\in X$. This defines a function, say $\mult_X:X \to (\mathbb{N}, \geq)$ which is upper semi-continuous (see e.g. \cite[Section 2)]{Abad}, or \cite[Th 6.12)]{multiplicity}. 
Therefore the level sets define a partition into locally closed sets, which 
we call the stratification of $X$ defined by the multiplicity. More precisely, a stratum 
will be an irreducible component of a level set. Given $x\in X$, $S_{mult_X, x}$ will denote the stratum through the point $x$, which we view simply as a set in the topological space of $X$.

Let $\max \mult_X(\in \mathbb{N})$ denote the highest multiplicity at points of $X$, and let $\Max \mult_X$ be the closed set of points of highest multiplicity. As $\mult_X$ is upper-semi-continuous, after a suitable restriction of $X$ to a neighborhood of any point $x$, we may assume that
$S_{mult_X, x}=\Max \mult_X$.

Here $X_{red}$ will denote the reduced scheme, which is also pure dimensional and excellent, and has the same topological space as $X$. 
 A particular feature of the multiplicity, to be used along this paper, is the compatibility of this stratification with reductions. Namely the stratification defined by the functions 
$\mult_X$ and $\mult_{X_{red}}$, on $X$ and $X_{red}$ respectively, are the same; or say
\begin{equation}\label{MULTRED}
S_{mult_X, x}=S_{mult_{X_{red}}, x}
\end{equation}
(see e.g. \cite{Abad}, 2), or \cite{multiplicity}, Th 6.14). 

The multiplicity at a point $x$ is one if and only if $x\in X$ is regular (\cite[Th 40.6, p. 157]{Nagata}); and it follows that the stratification defined by the multiplicity on a pure dimensional scheme is trivial (i.e., the function $\mult_X$ is 
constant) if and only if $X_{red}$ is regular.

We say that a subscheme $Y\subset X$ is equimultiple at a point $x\in Y$ if $Y$ is regular at $x$ and $Y\subset S_{mult, x}$. A well known property of the multiplicity, to be discussed below, says that if 
$X\leftarrow X_1$ denotes the blow up at $Y$ then $\mult_X(x)\geq \mult_{X_1}(x_1)$ for all 
$x_1\in X_1$ mapping to $x$. 
In particular, 
\begin{equation}\label{ecxxm}
\max \mult_X\geq \max \mult_{X_1}
\end{equation} if $Y\subset \Max \mult_X$ is regular. This blow up is said to simplify the multiplicity if the inequality is strict.
We now formulate the {\em Simplification Problem (for the multiplicity) in dimension $d$}, say SPM(d),  as follows: Set $X$ as above, and $d=\dim X$. If $\mult_X$ is not constant,  construct 
\begin{equation}\label{SPd} X\leftarrow X_1 \leftarrow \dots \leftarrow X_r,
\end{equation}
as a sequence of blow ups at regular centers as before, so that $\max \mult_X> \max \mult_{X_r}$. 

If $X$ is a variety and if $\max \mult_{X_r}=1$ one obtains a resolution of singularities of $X$.
But we don't know of the existence of a simplification if $X$ is a variety over a field of positive characteristic.
 We want to study conditions, at least locally at a point $x\in X$, so that a simplification of the multiplicity $X$ can be constructed if we assume the existence 
of a simplification for any scheme of dimension $d'$, for $d'$ strictly smaller then $d$. In other words we will give conditions that ensure that, at least locally, a simplification can be obtained by induction on the dimension $d$. 
To this end we first approximate $X$ by a complete intersection of the same dimension, and then we produce a scheme of dimension $d'(<d)$ which will ultimately lead us to the result.
\end{parrafo}
\begin{parrafo} {\em Reduction of the simplification of the multiplicity to complete intersections.}

Suppose that $X=\Spec(k[X_1, \dots,X_n]/\p)$,
where $\p$ is a prime of height $n-d$, one can choose 
\begin{equation}\label{approx1}\{ f_1, \dots, f_{n-d}\}\subset \p \mbox{ , or say}
\end{equation}
\begin{equation}\label{approx2} X \subset X'=\Spec(k[X_1, \dots,X_n]/\langle f_1, \dots, f_{n-d}\rangle),
\end{equation} and $X'$ is a complete intersection. 
We show that locally at $x\in X$ one can {\em approximate} $X$, by a suitable  complete intersection $X'$, 
so that the existence of a simplification of $X$ is equivalent 
to that of a simplification of $X'$. 

Here we start with $x\in X$, $X$ of dimension $d$, and locally at $x$ we construct a complete intersection $X'$  so that the following properties hold:
\begin{enumerate}
\item $X\subset X'$, both have the same dimension, and
\begin{equation}\label{forma0} 
\Max mult_{X'}=\Max mult_{X}.
\end{equation}
\item Any sequence of blow ups at equimultiple centers over $X'$ induces naturally a diagram 
\begin{equation}
\label{gdigraf}
\xymatrix@R=0pt@C=30pt{
X '    & \ar[l]
 X'_1  \ar[l]     &  \ar[l]  \ldots  &  \ar[l] X'_r   \\
\\
\\
\\
X \ar[uuuu]  & 
 \ar[l] X_1\ar[uuuu]     \ar[l]    & \ldots   \ar[l] &   \ar[l] X_r\ar[uuuu]
}
\end{equation}
where the lower row is a sequence of blow ups at equimultiple centers, and all vertical morphisms are closed immersions. 
\item If $\max \mult_{X'} =\dots = \max \mult_{X'_r}$, then $\max \mult_{X} =\dots = \max \mult_{X_r}$, and 
\begin{equation}\label{forma1} 
\Max mult_{X'_i}=\Max mult_{X_i} , \ 0\leq i \leq r.
\end{equation}
\end{enumerate}

This shows that a simplification of $X'$ induces one of $X$, and vice versa. In particular it shows that the simplification problem of the multiplicity reduces to the case of complete intersections. \end{parrafo}
\begin{parrafo}
Our next result is {\em inductive} and {\em local}. It will require a condition on the tangent cone of the point:
$C_{X,x}=\Spec(gr_{\calo_X }(\M_x)).$
More precisely on $(C_{X,x})_{red}$ (the reduced scheme of $C_{X,x}$).
The proof will rely on the previous approximation of $X$ by a complete intersection.
\begin{theorem}\label{T3}{\bf (Main Theorem)}  Let $X$ be a pure dimensional excellent scheme, of dimension $d$ at a point $x$. Assume that $\calo_{X,x}$ contains a field and that the residue field $k(x)$ is perfect. If $(C_{X,x})_{red}$ is not 
regular, then there is a simplification
of the multiplicity of $X$, locally at $x$, provided we assume the existences of simplification of multiplicity for schemes of dimension $d'<d$.
\end{theorem}

\end{parrafo}
\begin{parrafo}{\bf Multiplicity vs Hilbert Samuel.}
So far we have fixed a pure dimensional scheme $X$, and the function $\mult_X:X \to \mathbb{N}$
was defined by assigning to $x\in X$ the multiplicity of $\calo_{X,x}$.
One can also assign to this local ring another invariant, which is a sequence of positive integers, say $l_x: \mathbb N \to \mathbb N $, setting $l_x(n)=length(\calo_{X,x}/\M^{n+1}_x)$. This is known as the Hilbert-Samuel sequence at $x\in X$ (i.e., at $\calo_{X,x}$), and we define now the Hilbert-Samuel function on $X$, say $HS_X:X \to \mathbb{N}^\mathbb{N}$,  $HS_X(x)=l_x\in \mathbb{N}^\mathbb{N}$.
This function can be suitably modified, by adapting the value at each point $x\in X$ in terms of the local dimension, 
in such a way that this modified function, which we call again $HS_X:X \to \mathbb{N}^\mathbb{N}$, is upper-semi-continuous when $\mathbb{N}^\mathbb{N}$ is ordered lexicographically.
%
%
%

In particular a stratification on $X$ is obtained by considering the level sets of this function, called the Hilbert Samuel stratification of $X$. 
Let $S_{HS, x}$ denote the stratum containing $x\in X$, known as the Samuel Stratum through $x$. Samuel proved in \cite{Samuel} that $HS_X(x)=l_x\in \mathbb{N}^\mathbb{N}$ encodes also the multiplicity of $X$ at $x$.
This shows that $S_{HS, x} \subset S_{mult, x}$, so the stratification of the Hilbert-Samuel function is a refinement of that obtained from the multiplicity. 
When $X$ can be embedded as a hypersurface in a regular scheme $W$, the stratifications on $X$ defined 
by the function $HS_X:X \to \mathbb{N}^\mathbb{N}$ is the same as that defined 
by $\mult_X:X \to \mathbb{N}$. But in general the stratifications will be different.

Fix $x\in Y \subset X$ where $Y$ is irreducible and regular at $x$. Let $X\leftarrow X_1$ denote the blow up at $Y$, and fix $x_1\in X_1$ mapping to $x$.
Hironaka gives an algebraic characterization of the condition $Y\subset S_{HS, x}$ locally at $x$, know as the condition of {\em Normal Flatness of $X$ along $Y$ at the point $x$}. Namely, that the algebra $gr_X(I(Y))=\oplus I(Y)^n/I(Y)^{n+1}$ (the ring of functions on the normal bundle) be flat 
over $\calo_Y=\calo_X/I(Y)$, locally at the point $x$ (\cite{Hironaka64},Theorem 2, p. 195).

Theorems of Hironaka and Bennett (\cite{BB}), later simplified by Singh (\cite{BS}), say that under these conditions $HS_{X}(x) \geq HS_{X_1}(x_1)$. 
In particular, if $Y$ is regular and $Y\subset \Max HS_X$, then 
\begin{equation}\max HS_X\geq \max HS_{X_1}.
\end{equation}
This parallels the inequality in (\ref{ecxxm}). Therefore, given $X$ of dimension $d$, one can formulate the Simplification Problem in dimension $d$ for the Hilbert-Samuel function, in analogy with (\ref{SPd}). This was solved by Hironaka for schemes containing a field of characteristic zero. Let us simply indicate that his proof does not follow by induction on $d=\dim X$. Moreover, we do not know of a statement as that of the inductive Theorem \ref{T3} if {\em multiplicity} is replaced by {\em normal flatness}, not even in characteristic zero.
In fact Hironaka's proof used an embedding $X\subset W$ where $W$ is regular, and he argues by induction on the dimension of $W$ (not on the dimension of $X$, see \cite[p. 177]{Hironaka64}).
We refer to \cite{CJS} for a nice introduction to the Hilbert Samuel function. 
\begin{parrafo}{\em Some features of the multiplicity.}
If $X'\subset W$ is a hypersurface in a regular scheme, then, as mentioned, $\Max HS_{X'}=\Max \mult_{X'}$, and the same holds for the blow up, and for any sequence of blow ups. In fact for a hypersurface $X'$, the construction of a sequence of blow ups at normally flat centers is the same as one with equimultiple centers. 
This property will also hold for the complete intersection $X'$ in (\ref{forma0}), namely $\Max mult_{X'} =\Max HS_{X'}$. Another particular feature of the scheme $X'$ that we will construct will be that it {\em remains a complete intersection} after blowing up at any regular center included in $\Max mult_{X'} (=\Max HS_{X'})$. In summary, the complete intersection $X'$ assigned to $X$ will be very similar to a hypersurface. This will enable us to reformulate (\ref{forma1}):
\begin{equation}\label{forma10} 
\Max HS_{X'_i}=\Max mult_{X'_i}=\Max mult_{X_i} , \ 0\leq i \leq r.
\end{equation}

Given $X$, the existence of complete intersections $X'$ with the previous properties was presented  
in \cite{multiplicity}, where resolution of singularities (in characteristic zero) is proved by using the multiplicity instead of the Hilbert function (see also \cite{Abad}, \cite{BrV1}, and \cite{Nob}). The approximation of $X$ by $X'$ with the properties (\ref{forma0}), (\ref{gdigraf}), and (\ref{forma1}), to be discussed in Section 5, hold for the multiplicity (see  right hand term in (\ref{forma10})). Similarly for Theorem \ref{T3}, which also follows from this approximation.


Let us indicate here that historically a first and major step in the study of the multiplicity and its behavior under blow ups is due to Dade (\cite{D}): 
A subscheme $Y\subset S_{mult, x}$ which is regular at $x\in Y$ is said to be  {\em equimultiple along $Y$ at the point $x$}.  
He characterizes equimultiplicity by the condition that all fibers $\pi^{-1}(y)$ have same dimension, for $y\in Y$ in a neighborhood of $x$. 
He also proves that under such conditions $\mult_X(x) \geq \mult_{X_1}(x_1)$ at any $x_1$ mapping to $x$ (\ref{Dade}).

Later Hironaka and Schickhoff gave a second characterization of equimultiplicity in \cite{HironakaS} and \cite{Sch} respectively, which is also very geometric. To motivate the idea 
assume that $x\in Y \subset X \subset W$, where $Y$ and $W$ are regular, and that there is a {\em regular} subscheme, $H\subset W$, of complementary dimension with $Y$ and cutting $Y$ transversally at $x$. Consider the blow up, $W \stackrel{\pi'}{\longleftarrow} W_1$ of $W$ at $Y$, and let $H_1$, $X_1$ be the strict transforms, in $W_1$, of $H$ and $X$ respectively. So we obtain a diagram of blow ups, say
\begin{equation}\label{eq1}
\begin{array}{ccccccc}
 \overline{X}_1&\subset& H_1 & \subset  &   W_1&\supset &X_1\\
 \overline{\pi} \downarrow&&\downarrow  &   &  \downarrow \pi'& &\downarrow \pi\\
 \overline{X}=X\cap H&\subset& H& \subset  &  W & \supset& X
\end{array}
\end{equation}
By restriction of $\pi'$, we get $X \stackrel{\pi}{\longleftarrow} X_1$
 (the blow up of $X$ at $Y$) and also $H \leftarrow H_1$ (the quadratic transformation of $H$ at $x$). 
 Finally let
  $\overline{X}\stackrel{\overline{\pi}}{\longleftarrow}  \overline{X}_1(\subset H_1)$ be the quadratic transform of the section $\overline{X}$ at $x$. 
The diagram shows that there is a natural inclusion $\overline{X}_1\subset X_1\cap H_1$, and hence there is a natural inclusion  $\overline{\pi}^{-1}(x) \subset 
\pi^{-1}(x)$. Equimultiplicity of $X$ at $Y$ locally at $x$ is characterized by the condition $\overline{X}_1= X_1\cap H_1$ as sets, or equivalently, when  $$(\overline{\pi}^{-1}(x))_{red} =
(\pi^{-1}(x))_{red}$$ (see \ref{Par28}).
So this characterization expresses the equimultiplicity of $Y$ at $x$ in terms of the blow up of a transversal section $X\cap H$. Let us indicate that the inclusion in a regular scheme $W$ will be irrelevant, and also the choice of the section $H$. We include a proof of this result in Theorem \ref{TH211}, where we follow essentially that of Lipman in \cite[Section 5]{Lipman2}. The arguments we use for the proof will be crucial for our further discussion on the behavior of the multiplicity under blow ups, and for the proof of Theorem \ref{T3}.
\end{parrafo}
{\bf Multiplicity vs Hilbert Samuel at the tangent cone.}
There are many local invariants at a point $x\in X$ that can be reformulated in terms of the origin of the tangent cone, say $\mathbb{O}\in C_{X,x}$. For examples the multiplicity and the Hilbert-Samuel function: Both invariants give the same information 
for $\calo_{X,x}$ and for $\calo_{C_{X,x}, \mathbb{O}}$.

In this paper a cone over a field $k$, say $C=\Spec(\mathcal{A})$, is the spectrum of a graded $k$-algebra $\mathcal{A}$ which is generated in degree one. 
Hironaka studied 
the Hilbert Samuel stratum of a cone 
at the origin $\mathbb{O}\in C$, say $S_{HS, \mathbb{O}}$. He proved two results
\begin{enumerate}
\item that every time we fix an embedding (of cones), say $C\subset \mathbb{V}$ in a vector space over $k$, the stratum $S_{HS, \mathbb{O}}=\mathbb{S}$ is a subspace of $\mathbb{V}$,
\item this subspace acts on the cone, namely $C_{}+\mathbb{S}=C_{}$ in $ \mathbb{V}$.
\end{enumerate}
To be precise, $\mathbb{S}$ is a subspace when  the characteristic of $k$ is zero. However this property has a natural extension to positive characteristic. This extension requires the study of other subgroups in $ \mathbb{V}$, which are not necessarily subspaces. This led Hironaka to the use of 
group schemes in positive characteristic, also studied in work of Giraud, Oda, and Pomerol among others (see \cite{Giraud1975}, \cite{GiraudNotas}, \cite{Hironaka70}, \cite{Oda1987}, \cite{Oda2}, \cite{Pomerol}). All these works are directed to the study of the Hilbert-Samuel functions along a singular variety, whereas here 
we will consider the multiplicity as main invariant. 
Concerning the study of normal flatness in positive characteristic there is also a recent work of Dietel in \cite{Dietel}.

Summarizing, groups schemes turn out to be the appropriate language to express some notions as {\em vector spaces}, and the action of a subspace.
If $C_{X,x}\subset T_{X,x}$ is the inclusion of the tangent cone in the Zariski tangent space, group schemes enables us to endow the scheme $T_{X,x}$ with a natural structure of vector space.
We shall briefly recall results of this theory, at least those to be used in this paper, and we  give precise reference for the reader not acquainted with this concept.

We will show here that if $C$ is a cone over a perfect field (i.e., if the residue field at $\mathbb{O}\in C$  is perfect), then the stratum of the multiplicity $S_{mult, \mathbb{O}}(\subset C)$  is also the Hilbert-Samuel stratum of the reduced scheme $C_{red}$ through the origin (see \ref{lem51}).

The following properties of the multiplicity, gathered in the next theorem,  are related to theorems of Hironaka concerning the notion of normal flatness. In other words, Hironaka studies similar or related properties when the invariant is the Hilbert Samuel invariant, and when the blow up is defined at normally flat center $Y$: 
%
%

\begin{theorem}\label{T1} Fix a point $x\in X$ and let $C_{X,x}\subset T_{X,x}$ be the inclusion of the tangent cone in the Zariski tangent space. Assume that the residue field $k(x)$ is perfect, and let $\mathbb{S}$ denote the stratum of (highest) multiplicity through the origin of the cone 
$C_{X,x}$.
\begin{enumerate}
\item $\mathbb{S}$ is a subspace in $T_{X,x}$. It acts on $(C_{X,x})_{red}$, and it is the biggest subspace with this property. 

\item Let $Y\subset X$ be regular and equimultiple at $x$, then the subspace 
$T_{Y,x}$ is included in $\mathbb{S}$, and hence it also acts on the subscheme $(C_{X,x})_{red}$. In addition, this action provides a natural identification $(C_{X,x})_{red}/T_{Y,x}=(C_{X,Y,x})_{red},$ where $C_{X,Y,x}$ is the normal cone of $x\in Y\subset X$.
\item Let $X\stackrel{\pi_Y}{\longleftarrow} X'$ be the blow up at $Y$, as above. Given $x'\in \pi^{-1}(x) \subset X'$, 
$e(\calo_{X', x'}) \leq e(\calo_{X, x})$, and if the equality holds then
\begin{equation}x' \in \Proj (\mathbb{S}/T_{Y,x})\subset \pi^{-1}(x)=\Proj(C_{X,Y,x}),
\end{equation}
where the inclusion is that derived from the action of $T_{Y,x}$ on $\mathbb{S}\subset (C_{X,x})_{red}$ in $T_{X,x}$.
%
\end{enumerate}
\end{theorem}
The statement in (1) does not hold if $k(x)$ is not perfect (\cite[Ex 2.12 , III-25]{GiraudNotas}) . The inequality $e(\calo_{X', x'}) \leq e(\calo_{X, x})$ in (2) is Dade's result mentioned before (\ref{Dade}).

\vskip 3cm
So according to this result, if $dim(\mathbb{S})=0$, setting $X\leftarrow X_1$ is the quadratic transform at $x$, then, after a suitable restriction of $X$ to an open neighborhood of $x$, $\max \mult_X> \max \mult_{X_1}$.

Furthermore, there is information extracted from $\mathbb{S}$ that will rule the behavior of the singularity at $x\in X$, when applying {\em any sequence} of blow ups at regular {\em equimultiple} centers:
\end{parrafo}
\begin{theorem}\label{T2} Assume that the residue field $k(x)$ is perfect at a point $x\in X$ and let $\mathbb{S}$ denote the stratum of (highest) multiplicity through the origin of the cone 
$C_{X,x}$. Assume, for simplicity that 
$x\in \Max mult_{X} $, and let $D(x)$ be the dimension of the subspace $\mathbb{S}\subset T_{X,x}$
\begin{enumerate}
\item $D(x)$ is an upper bound 
of the local dimension of the closed set $\Max mult_{X} $ at $x$.

\item For any sequence $X\leftarrow X_1 \leftarrow \dots \leftarrow X_r$, of blow ups at regular equimultiple centers, and given points $x_i\in X_i$, $0\leq i \leq r$ so that $x_{i+1}$ maps to $x_i$, and 
$x_0=x$, if $e(\calo_{X_r, x_r})=e(\calo_{X, x})$, then $D(x)$ is also an upper bound 
of the local dimension of the closed set $\Max mult_{X_r} $ at $x_r$.
\end{enumerate}
\end{theorem}

So for example, if $D(x)=1$, then any sequence over $X$, defined as in the Theorem, consists on blowing up
either points or regular curves.

Our Main Theorem \ref{T3} will follow from Theorems \ref{T1} and \ref{T2}.
The statements in \ref{T2} parallel a Theorem of Hironaka stated here as Theorem \ref{T22}. However in Theorem \ref{T22} the centers are chosen to be normally flat, whereas here in \ref{T2}
centers are equimultiple.
\begin{parrafo} The paper is organized as follows: In Sections 2 and 3 we prove the two characterizations of equimultiplicity mentioned before.
The techniques introduced there will be used along the paper.
In Section 4 we discuss the stability of transversality when blowing up at equimultiple centers. In Section 5 we show that a scheme can be replaced by a nice complete intersection if we try to find a simplification of the multiplicity.

The notion of group schemes is briefly discussed in Section 6, and applied in the proof of 
Theorem \ref{TH57} which concerns the stratification defined by the multiplicity on any affine cone.
Finally, Theorems \ref{T1} and \ref{T2} are proved in \ref{PT1} and \ref{PT2} respectively, and the Main Theorem in \ref{PT3}.

As indicated above, Theorems \ref{T1} and \ref{T2} 
extend to the multiplicity results which are already known 
for the Hilbert-Samuel function. Whereas the inductive result in Main Theorem  \ref{PT3} is based on the {\em approximation} by complete intersections, introduced in \cite{multiplicity}, which is exclusive of the multiplicity.

I profited from discussions with C. Abad, A. Benito, A. Bravo, S. Encinas, and D. Sulca.

\end{parrafo}
\end{section}
\begin{section}{Conditions for equimultiplicity and a theorem of Dade.}
In this section we recall some basic notions on integral closure of rings and ideals. The main 
result in this section is a theorem of Dade, in Theorem \ref{TH14}, which characterizes equimultiplicity
in terms of integral closure of ideals and the dimension of the fibers of the blow up.
%
\begin{parrafo}\label{finite}
Northcott and Rees introduce the notion of reduction. Given ideals $I\subset J$ in a noetherian ring $B$, $I$ is said to be a reduction of $J$ if $IJ^n=J^{n+1}$ for some integer $n$. Equivalently, if the inclusion of Rees rings $B[IW]\subset B[JW]$ is a finite extension of subrings in $B[W]$.

In this case, if $X_I\to \Spec(B)$ denotes the blow up at $I$ and $X_J \to \Spec(B)$ is the blow up at $J$, there is a factorization $X_J\to X_I$ which is induced by this finite extension, and hence it is also a finite morphism.

The notion of reduction of an ideal $J$ in $B$ will appear naturally when studying the fibers of the blow up
$X_J \to \Spec(B)$. Note that if $(A,M)$ is the localization of $B$ at a prime ideal, the fiber over the closed point is the projective scheme attached to $A/M\otimes_AA[JW]=A/M\otimes_A gr_A(J)$.
The following result give a useful criterion to produce a reduction, at least for an ideal in a local ring.

\begin{theorem}(\cite[Th.10.14]{Herrmann}\label{1.2} Let $(A,M)$ be a local ring with residue field $k$, and let $J$ be a proper ideal. For a given element $a \in J$, let $a^{*}$ denote the class of $a$ in $J\otimes_Ak=J/MJ$.
Given $a_1, \dots, a_s$ in $J$, the following conditions are equivalent:

i) $I=\langle a_1, \dots, a_s \rangle$ generate a reduction of $J$.

ii) The graded $k$-algebra $(k\otimes_A gr_A(J))/\langle a_1^{*}, \dots, a_s^{*} \rangle$ is zero dimensional.
\end{theorem}
\proof

Set $G=(k\otimes_A gr_A(J))/\langle a_1^{*}, \dots, a_s^{*} \rangle= (k\otimes_A (A[JT]))/\langle a_1^{*}T, \dots, a_s^{*} T\rangle$. The condition in (ii) is equivalent to $G_n=0$ for $n>>0$, where $G_n$ denotes the homogeneous component in degree $n$. 

If (i) holds, then $IJ^n=J^{n+1}$ for $n$ big, so (ii) holds. Conversely, if (ii) is satisfied, $ J^{n+1} = IJ^n+MJ^{n+1}$, so $IJ^n=J^{n+1}$, for $n$ big.
\end{parrafo}

\begin{parrafo}\label{finito}
When $k=A/M$ is infinite, one concludes from Noether's Theorem that for an integer, say
$e\geq  dim(k\otimes_A gr_A(J)))$, one can find $e$ elements in $J$ which span a reduction of the ideal. This requirement on the residue field is not a restriction for our purpose. In fact, the properties to be studied here, such as the multiplicity, are compatible with \'etale 
topology. So, when the residue field of $(A,M)$ is finite, one can argue at the strict henselization (which does have an infinite residue field), and finally descent to a suitable \'etale neighborhood.

A local ring $(A,M)$ is said to be formally equidimensional (quasi-unmixed in Nagata's terminology) if 
$\dim(\hat{A}/p)=\dim(\hat{A})$
at each minimal prime ideal $p$ in the completion $\hat{A}$.

A first connection of integral closure with the notion of multiplicity is given by the following Theorem of Rees.

\begin{theorem}\label{Rees} \cite{Rees} If $I\subset J$ are primary ideals for the maximal ideal in a formally equidimensional local ring $(A,M)$, then both ideals have the same integral closure if and only if $e_A(I)=e_A(J)$.
\end{theorem}
\end{parrafo}
\vskip 0,5cm

\begin{parrafo}(On multiplicity and finite extensions) The following theorem, crucial in our discussion, relates the behavior of the multiplicity of ideals under finite ring  extensions.
\begin{theorem}\label{MultForm} \cite[Theorem 24 p. 297]{ZS} Let 
$(A,M)$ be a local domain, and let $B$ be a finite extension of $A$. Let $K$ denote the quotient field of $A$, and $L=K\otimes_AB$.

 Let $Q_1, \ldots, Q_r$ denote the maximal ideals of the semi-local ring $B$, and assume that $\dim B_{Q_i}=\dim A$, $ i=1, \dots ,r$. Then
$$e_A(M)[L:K] =\sum_{1\leq i \leq r} e_{B_{Q_i}}(M B_{Q_i}) [k_i:k],$$
where $k_i$ is the residue field of $ B_{Q_i}$, $k$ is the residue field of $(A,M)$, and $[L:K]=\mbox{dim}_KL$.

\end{theorem}
\begin{parrafo}\label{CincluidoB}
We say that a ring is {\em pure dimensional or equidimensional} when every saturated chain of prime ideals has the same length. All rings to be considered here will be pure dimensional  and excellent, so the localization at any prime is formally equidimensional. Typically we will consider pure dimensional $k$-algebras over a perfect field, their localizations at prime ideals, their completions or their henselizations. 

In our discussion we consider a ring $B$, as above, and we assume that there is a regular subring $S\subset B$ so that the extension is finite. Let $K$ be the quotient field of $S$, and let $L=B\otimes_S K$. 
For example, if $X$ is of finite type over a perfect field $k$ and $x\in X$ is a closed point of multiplicity, say $n$, then after replacing the point by an \'etale neighborhood, one may assume that $X=\Spec(B)$, and there is a regular subring $S$ with 
the additional condition that $L=B\otimes_S K$ is of dimension $n$ over $K$ (see 
\cite[Appendix 1)]{BrV1}). We say that 
$B$ has generic rank $n$ over $S$. Similar statements hold for complete pure dimensional local rings containing a field.

Let $P$ be a prime ideal in $B$ and $\p=P\cap S$.  Let $P=P_1, P_2\ldots, P_r$ denote the prime ideals in $B$ that dominate $\p$ in $S$. Here we assume that $\dim B_{Q_i}=\dim S_{\p}$, $ i=1, \dots ,r$. Then, the previous Theorem together with usual properties of the multiplicity of ideals, show that
$$[L:K] =\sum_{1\leq i \leq r} e_{B_{P_i}}(\p B_{P_i}) [k_i:k] \geq e_{B_{P}}(\p B_{P})\geq e_{B_{P}}(P B_{P}).$$
\end{parrafo}

\begin{corollary}\label{C16} (of Th \ref{MultForm}) Let $S\subset B$ be a finite extension, where $S$ is a regular domain. Let $P$ be a prime ideal in $B$ and $\p=P\cap S$. If $B_P$ is formally equidimensional the following conditions 1) and 2) are equivalent:

1) $e_{B_P}(PB_P)= [L:K]$;

2) 

2i) $P$ is the only prime in $B$ dominating ${\p}$ (i.e., $B_P=B\otimes_SS_{\p}$);

2ii) $S_{\p}/{\p}S_{\p}=B_P/PB_P$;

2iii) ${\p}B_P$ is a reduction of $PB_P$ in $B_P$.
\end{corollary}
 The claim follows from Theorems \ref{MultForm}  and \ref{Rees}.

\begin{theorem}\label{T18} Given $S\subset B$ as before, and assume that $B_P$ is formally equidimensional at any prime $P$. Set $n=[L,K]$ and consider the finite morphism $\delta: \Spec(B)\to \Spec(S)$ and let $F_n(B)$ denote the set of primes $P$ where $B_P$ has multiplicity $n$. Then
\begin{enumerate}
\item $e_{B_P}(PB_P)\leq n$ at any prime $P$ of $B$.
\item The function $\delta: F_n(B) \to \delta(F_n(B))$ is a bijection of sets. In particular,
if $P\in F_n(B)$, $B_P=B\otimes_SS_{\p}$ where $\p=P\cap S$, and the conditions in 2) of the previous Corollary hold.
\item If $P\in F_n(B)$ and $B/P$ is regular, then  
$S/(P\cap S)=B/P$.
\end{enumerate}
\end{theorem}
The first two claims follow from the previous discussion. The third is a consequence of (ii) in the next Lemma.
\begin{lemma}\label{R16}  Let $Q\subset P$ be an inclusion of primes in $B$. Set $ \p=P\cap S$, $\overline{B}=B/Q$, $\overline{S}=S/(Q\cap S)$, 
$\overline{P}=P/Q$, and $\overline{\p}=\p / (Q\cap S)$. Then,
\begin{enumerate}
\item[i)] If the equivalent conditions in the previous Corollary \ref{C16} hold for $P$ and $\p$, they also hold for
$\overline{P}$ and $\overline{\p}$.
\item[ii)] If $B/Q$ is a regular ring, then $\overline{B}_{\overline{P}}=\overline{S}_{\overline{\p}}$.
\end{enumerate}
\end{lemma}
\proof The claim (i) is a consequence of the three conditions in 2) of the corollary. For (ii) note that
$\overline{S}_{\overline{\p}} \subset \overline{B}_{\overline{P}}$ a finite extension of domains, and 
$\overline{\p}\overline{B}_{\overline{P}}$ is a reduction of the maximal ideal of 
$\overline{B}_{\overline{P}}$. By assumption the local ring $\overline{B}_{\overline{P}}$ is regular, so the maximal ideal does not admit any proper reduction. Hence
$\overline{\p}\overline{B}_{\overline{P}}=\overline{P} \ \overline{B}_{\overline{P}}$. This, together with the other conditions in part 2) of the corollary show that 
$\overline{B}_{\overline{P}}$ has rank one over $\overline{S}_{\overline{\p}}$. This proves the equality (ii), so both rings are regular.
\end{parrafo}

\begin{parrafo} \label{110}{\em A generalization of Rees Theorem.}
Let $I$ denote an ideal in a local ring $(A,M)$. Let 
$f: X \to \Spec(A)$
be the blow-up at $I$, and let
$f_0: X_0\to \Spec(A/I)$
be the proper morphism induced by restriction.
Northcott and Rees defined the {\em analytic spread} of $I$ as: 
$$l(I)=dim(A/M\otimes_Agr_A(I))=\delta +1,$$ where $\delta$ is the dimension of the fiber of $f$ over the closed point of $\Spec(A)$, or equivalently, the dimension of the fiber of $f_0$ over the closed point.
Note that $l(I)=dim(A)$ when $I$ is $M$-primary.

The height of $I$, say $h(I)$, is  $min(\dim A_p)$ as $p$ runs through all primes containing $I$, and 
$$l(I) \geq h(I)$$
with equality if and only if all fibers of $f_0$ have the same dimension.
The inequality holds because the the dimension of the fibers of $f_0: X_0\to \Spec(A/I)$ is an upper semicontinuous function on primes of $A/I$. In addition, if $\p$ is minimal containing $I$, the dimension of  $f_0^{-1}(\p)$ is dim $A_{\p}$.

Let  $I\subset J$ be a reduction of ideals in a noetherian ring $B$, and  let $X_I\to \Spec(B)$ and $X_J \to \Spec(B)$ denote the blow ups at $J$ and $I$. Since there is a factorization $X_J\to X_I$ which is finite, it follows that $l(IB_P)=l(JB_P)$ at any prime $P$ in $B$.

We refer to Theorems 2 and 3 in \cite{Lipman2} for the following result of Böger.

\begin{theorem}\label{TB}  Let $I\subset J \subset \sqrt{I}$ be ideals in a formally 
equidimensional local ring $A$. If $h(I)=l(I)$, then $I$ is a reduction of $J$ if and only if 
 $$e_{A_{\p}}(IA_{\p})=e_{A_{\p}}(JA_{\p})$$
at every minimal prime $\p$ of $I$.
\end{theorem}

\end{parrafo}
\begin{parrafo}
We will draw special attention to the blow up of schemes at regular centers. Over fields of characteristic zero resolution of singularities can be achieved by blowing up at regular equimultiple centers. We first introduce a characterization of equimultiplicity at the completion of a local ring $(A,M)$, to be generalized later as a condition on the ring itself. Here equimultiplicity of a prime $P$ at $A$ means that $(A,M)$ and $A_P$ have the same multiplicity. One consequence that we will extract from the proof is that 
such $P$ will also fulfill the condition $ht(P)=l(P)$ (Theorem \ref{TH14}).
\color{black}
\end{parrafo}\color{black}
\begin{lemma}\label{L13} Let $(A,M,k)$ be an excellent pure dimensional local ring of dimension $d$, containing a field, and assume that the residue field $k$ is infinite. Let $P$ be a prime such that $(A/P,M/P)$ is regular. Let $(B,M')$ be the completion of $(A,M)$, and $P'=PB$. The following conditions are equivalent:  

i) $e=e_A(M)=e_{A_P}(PA_P)$.

ii) There is a family $y_1, \dots, y_d$ in $B$  and an integer $1\leq s \leq d$, so that 
$\langle y_1, \dots, y_d \rangle B$ is a reduction of $M'$ and $\langle y_1, \dots, y_s \rangle B$
is a reduction of $P'=PB$.

In addition, if ii) holds the generic rank of $S=k[[y_1, \dots, y_d]]\subset B$ is the multiplicity $e$.

\end{lemma}
\proof The graded ring $gr_A(M)$ is a $k$-algebra of dimension $d$, and as $k$ is  infinite, Noether's normalization ensures that there is an inclusion of graded rings
$k[X_1, \dots, X_d] \subset gr_A(M)$ which is a finite extension, for a suitable choice of $X_i\in M/M^2$. Let $x_i\in M$ be an element with class $X_i$ at $M/M^2$. 
By construction, and by Theorem \ref{1.2},
$\langle x_1, \dots, x_d\rangle$ is a reduction of $M$ in $A$, and there is an inclusion 
$S=k[[x_1, \dots, x_d]]\subset B,$
where $B$ is the completion of $A$, and $S\subset B$ is finite.
Let $K \subset L$ denote the quotient fields of $S$ and $B$. 
Excellence ensures that $e_A(M)=e_B(M')$. Theorems \ref{Rees} and \ref{MultForm} show that $[L,K]=e=e_A(M)$; namely, that the generic rank $[L:K]$ is the multiplicity of $(A,M)$  (or say of the completion $(B,M'))$, and hence the multiplicity at any prime ideal of $B$ is at most $e$. If  $P$ is a prime of multiplicity $e$, namely if $A_P$ has multiplicity $e$, and $P$ is a regular prime in $A$ (i.e., if $A/P$ is regular), it induces a regular prime, say $P'$ in $B$, and Lemma \ref{R16} says that $\p=P'\cap S$ is also regular.

Let $\{x'_1, \dots, x'_d\}$ be a regular system of parameters of $S$ such that $\p=\langle x'_1, \dots, x'_h \rangle$ for some $h\leq d$.

If i) holds, Theorem \ref{T18} ensures that $P'$ is the unique prime in $B$ dominating $\p$. This can be checked replacing $S$ by $S_{\p}$. In addition, $\p B_{P'}$ is a reduction of $P'B_{P'}$ (see Corollary \ref{C16}).

We claim now that $\p B=\langle x'_1, \dots, x'_h \rangle B$ is a reduction of $P'$ in $B$.
Note first that $ht(\p B)=l(\p B)$, in fact the blow-up of $B$ at $\p B$, is finite over the blow up of the regular ring $S$ at the regular prime $\p$. In particular the fibers of the two blow ups have the same dimension; and clearly $ht(\p S)=l(\p S)$. So the claim follow from Theorem \ref{TB}.

Finally, let $x'_{h+1}, \dots, x'_{d}$ be elements in $A$ inducing a regular system of parameters in $A/P$, and note that $\langle x'_1, \dots, x'_h, x'_{h+1}, \dots , x'_d \rangle B$ is a reduction of $M'$, as claimed in (ii).

For the converse, if (ii) holds consider the finite extension $S'=k[[y_1, \dots, y_d]]\subset B.$
As $\langle y_1, \dots, y_d \rangle $ is a reduction of $M'$, we conclude that $[L: K']=e$, where $K'$ denotes the quotient field of $S'$.
As $\langle y_1, \dots, y_h \rangle B$ is a reduction of $P'$ in $B$, it is the only prime in $B$ dominating $S'$ at $\p'=\langle y_1, \dots, y_h \rangle S'$. Finally, the conditions in 2) of Corollary \ref{C16} hold, and hence $e_{B_{P'}}(P'B_{P'})=[L, K']=e$.

Excellence ensures that the morphism $\Spec(B)\to \Spec(A)$ is regular, in particular primes in correspondence have the same multiplicity. So $e_{A_{P}}(PA_{P})=e_{B_{P'}}(P'B_{P'})=e=e_{A}(M)$.
\endproof

We now present a proof of a result of Dade (see also \cite[Corollary p.121]{Lipman2}).

\begin{theorem}\label{TH14} Let $(A,M)$ be an excellent pure dimensional local ring of dimension $d$, containing a field, and assume that the residue field is infinite. Let $P$ be a prime such that $(A/P,M/P)$ is regular. The following conditions are equivalent. 

i) $e=e_A(M)=e_{A_P}(PA_P)$;

ii) $l_{(A,M)} (P)=ht(P)$.
\end{theorem}

\proof i) implies ii) is a consequence of the previous lemma. There we shows that if (i) holds, there are elements $y_1, \dots, y_d$ in $B=\hat{A}$  and an integer $1\leq s \leq d$, so that 
$\langle y_1, \dots, y_d \rangle$ is a reduction of $M'$ and $\langle y_1, \dots, y_s \rangle$
is a reduction of $P'=PB$. 
Set $S=k[[y_1, \dots, y_d]]\subset B$ as above, and note that the blow up of $B$ at the ideal $\langle y_1, \dots, y_s \rangle B$ is finite over the blow-up of $S$ at the regular prime $\langle y_1, \dots, y_s \rangle$; so the closed fibers of both blow ups have the same dimension, and the equality in ii) is now clear (see \ref{110}).

Let $s=ht(P)$ and assume that (ii) holds (i.e., that  the dimension of  
$A/M\otimes_A gr_A(P)$ is $s$). Theorem \ref{1.2} states that one can choose 
$s$ elements $y_1, \dots, y_s$ in $P$, which span a reduction of $P$. By assumption 
$dim(A/P)+dim(A_P)= d$. Finally, since $A/P$ is regular, one can extend $y_1, \dots, y_s$, to say $y_1, \dots, y_s, \dots, y_d$ which clearly span 
a reduction of $M$, so (i) follows from Lemma \ref{L13}.
\endproof

\end{section}
\begin{section}{Equimultiplicity and a theorem of Hironaka-Schickhoff}
In this section the main result is the characterization of equimultiplicity in Theorem \ref{TH211}. This characterization, discussed in the Introduction (see \ref{eq1}), is very geometrical. Here we follow the algebraic proof due to J. Lipman, with some minor changes.  Some previous technical  results will be needed, also for our further discussions, and we suggest here a first look at \ref{Par28}, \ref{Par38}, and to the formulation of the theorem.

\begin{remark} The integral closure of ideals in a ring $B$ is an operation which is compatible with inclusions of ideals. The integral closure of zero ideal are the nilpotent elements, and an element $\theta$ is in the integral closure of an ideal $J$ if and only if the same holds at the reduced ring $B_{red}$. 

We will assume here that $B$ an excellent pure dimensional {\em reduced} ring  and  $\widetilde{B}$ will denote its integral closure. Here $\widetilde{J}$ will denote the integral closure in the ring $B$, and $\widetilde{J}'$ will denote the integral closure in $\widetilde{B}$ of the extended ideal $J\widetilde{B}$, so $\widetilde{J}=B\cap \widetilde{J}'$.
The inclusion $$B \oplus J \oplus J^2 \oplus \dots \subset 
\widetilde{B} \oplus \widetilde{J}'\oplus \widetilde{J^2}' \oplus \dots$$
is a finite extension, and furthermore it is a normalization. Hence the projective scheme defined by this latter ring is the normalization of the blow up defined by the former.

Let $X'\to X$ be any proper dominant morphism of reduced schemes, and let $I$ be an ideal in 
$\calo_X$. Let $\widetilde{I^{e}}$ denote the integral closure of the extended ideal $I^{e}=I\calo_{X'}$. Note that  $\widetilde{I^{e}} \cap \calo_X$ is the integral closure of $I$ in $\calo_X$.
This claim follows from the characterization of the integral closure of ideals, and also of proper morphisms, in terms of valuation rings.

\end{remark}

\begin{parrafo}\label{impcri} We consider now rings and ideals in the following situation:
\begin{itemize}
\item $S\subset B$ is a finite extension of reduced pure dimensonal excellent rings (not necessarily local).
\item $S$ is a regular domain and $N$ is a regular prime in $S$ (i.e., $S/N$ is regular).
\item $Q$ will denote the extended ideal $NB$.
\end{itemize}
Note that every non-zero element in $S$ is a non-zero divisor in $B$. If $K$ denotes the quotient field $S$, then $B$ is included in $ L=B\otimes_SK$, which is the total quotient ring, and the minimal polynomial over the field $K$ of an element $\theta \in B$ has coefficients in $S$ since the latter ring is normal. The following will be a criterion to in our discussion.
\end{parrafo}
\begin{proposition}\label{prop45} 
Fix $\theta \in B=B_{red}$,  with minimal polynomial $V^m+a_1V^{m-1}+ \dots+a_m\in S[V]$. Then $\theta \in \widetilde{Q^r}$ (the integral closure of $Q^r$ in $B$) if and only if 
\begin{equation}\label{eq44} \frac{\nu_N(a_i)}{i}\geq r
\end{equation}
for $i=1, \dots, m.$
\end{proposition}
%
\proof As $N$ is regular in $S$, powers and symbolic powers of $N$ coincide. So the conditions in (\ref{eq44}) ensure that $a_i\in N^{ri}$, and hence $\theta 
\in \widetilde{Q^r}$. 
For the converse consider the blow up
$$ \Spec(B) \longleftarrow X=\Proj (B \oplus Q \oplus Q^2 \oplus \dots)$$
and let $\mathcal{L}$ denote the invertible ideal $Q\calo_X=N\calo_X$.
The assumptions are such that 
$$S \oplus N \oplus N^2 \oplus \dots \subset B \oplus Q \oplus Q^2 \oplus \dots \subset \widetilde{B} \oplus \widetilde{Q}' \oplus \widetilde{Q^2}'
\oplus \dots$$
are finite extensions, where the latter is the normalization of the middle ring. 
Set 
$Z=\Proj(S \oplus N \oplus N^2 \oplus \dots )$, and let $\mathcal{L}_Z$ denote the invertible ideal $N\calo_Z$.
Set $\widetilde{X}=\Proj( \widetilde{B} \oplus \widetilde{Q}'\oplus \widetilde{Q^2}'
\oplus \dots)$. So there are finite morphisms $\widetilde{X} \to X\to Z$, and 
$\mathcal{L}_Z\calo_X=\mathcal{L}$ (see \cite[4.3, p. 336]{multiplicity}). 

Clearly $Q\widetilde{B}$ and the extension of $\widetilde{Q}$ to $\widetilde{B}$ have the same integral closure $\widetilde{Q}'$ in $\widetilde{B}$. As both are invertible at the normal scheme $\widetilde{X}$, also
$$\mathcal{L}_Z\calo_{\widetilde{X}}=\widetilde{Q}\calo_{\widetilde{X}}.$$ 
In fact, an invertible ideal in an integrally closed scheme is integrally closed.

If we replace $N\subset Q \subset \widetilde{Q}'$ (in $S\subset B \subset \widetilde{B}$) by
$N^r\subset Q^r \subset \widetilde{Q^r}'$ in the previous argument, then
$$\mathcal{L}^r_Z\calo_{\widetilde{X}}=\widetilde{Q^r}\calo_{\widetilde{X}}.$$

So if $\theta \in \widetilde{Q^r}(\subset B)$, we note that $\theta$ is a global section of $(\mathcal{L}_Z)^r\calo_{\widetilde{X}}$. 

As $\widetilde{X} \to X\to Z$ are affine, an affine cover of $Z$ induces one in $X$ and also in $\widetilde{X}$. One can fix finite extensions of rings $S_1\subset B_1 \subset \widetilde{B}_1$ so that 
$\Spec(S_1)$, $\Spec(B_1)$, and $\Spec(\widetilde{B}_1)$ are affine charts of $Z$, $X$, and $\widetilde{X}$ respectively, and the restriction of $\mathcal{L}_Z$ to $S_1$ is a principal ideal, say $x_1S_1$. As $\widetilde{B}_1$ is normal, 
$x_1\widetilde{B}_1$ is integrally closed in $\widetilde{B}_1$ and 
$\theta=x_1^r\cdot (\theta')$
for some element $\theta' \in \widetilde{B}_1$. In particular, as $\widetilde{B}_1$ is finite over $S_1$,
both $\theta$ and $\theta'$ are integral over $S_1$. Note that the minimal polynomial of  $\theta'$ is 
$$V^m+\frac{a_1}{x_1^r}V^{m-1}+ \dots+\frac{a_m}{x_1^{rm}}\in K[V],$$
where $K$ denotes the quotient field of $S$.
Since $S_1$ is regular, and therefore normal, this equation is also in $S_1[V]$. In particular 
it follows that $\nu_N(a_i)\geq mi$, for $i=1, \dots, m$.
\endproof

\begin{corollary}\label{cor47} Assume now that in the previous setting $S \subset B$ is an inclusion of local rings. Let $\{ x_1, \dots , x_d\}$ be a regular system of parameters in $S$. Fix integers 
$ e_1, e_2$ so that $1\leq e_1 < e_2 \leq d$. Let $Q_1, Q_2$ be two ideals in $B$, where 
$Q_1$ is integral over $\langle x_1, \dots, x_{e_1} \rangle B$, and $Q_2$ is integral over $\langle x_{e_2}, \dots, x_{d} \rangle B$. Then
$\widetilde{Q_1Q_2}=\widetilde{Q_1} \cap \widetilde{Q_2}$
(in particular $Q_1\cap Q_2 \subset \widetilde{Q_1Q_2}$).
\end{corollary}
\proof The inclusion $\widetilde{Q_1Q_2}\subset \widetilde{Q_1} \cap \widetilde{Q_2}$ is clear.
Fix $\theta \in \widetilde{Q_1} \cap \widetilde{Q_2}$  
with minimal polynomial $V^m+a_1V^{m-1}+ \dots+a_m\in S[V]$. Then, applying the proposition for $r=1$:
$$ a_i\in (\langle x_1, \dots, x_{e_1} \rangle)^i \cap (\langle x_{e_2}, \dots, x_{d} \rangle)^i=
(\langle x_1, \dots, x_{e_1} \rangle)^i \cdot (\langle x_{e_2}, \dots, x_{d} \rangle)^i=
(\langle x_1, \dots, x_{e_1} \rangle \cdot \langle x_{e_2}, \dots, x_{d} \rangle)^i.$$

\begin{remark}\label{rk26} 1) The previous conclusion applies also for powers of the two ideals. Namely, given positive integers $r$ and $s$, 
$\widetilde{Q^r_1Q^s_2}=\widetilde{Q^r_1} \cap \widetilde{Q^s_2}$ (in particular $Q^r_1\cap Q^s_2\subset \widetilde{Q^r_1Q^s_2}$). 
In fact, note that $$ (\langle x_1, \dots, x_{e_1} \rangle)^{ri} \cap (\langle x_{e_2}, \dots, x_{d} \rangle)^{si}=
(\langle x_1, \dots, x_{e_1} \rangle^r)^i \cdot (\langle x_{e_2}, \dots, x_{d} \rangle^s)^i=
(\langle x_1, \dots, x_{e_1} \rangle^r \cdot \langle x_{e_2}, \dots, x_{d} \rangle^s)^i.$$

2) Note that the equalities in 1), namely $\widetilde{Q^r_1Q^s_2}=\widetilde{Q^r_1} \cap \widetilde{Q^s_2}$ in the ring $B$, also hold if we replace $Q_1$ (or $Q_2$)  by an ideal with the same integral closure.
\end{remark}

\begin{proposition}\label{P27} Let $Q_1, Q_2$ be two ideals in $B$ in the conditions of Corollary \ref{cor47}. Set 
$\overline{B}=B/Q_1$, and  $\overline{Q}_2=Q_2 B/Q_1$. 
Then: 
\begin{enumerate}
\item The kernel of the natural surjection 
$B/Q_1\otimes_B gr_B(Q_2) \to gr_{\overline{B}}(\overline{Q}_2)\to 0$
is nilpotent.
\item The previous statement also holds if $Q_1$ or $Q_2$ are replaced by their integral closure.
\end{enumerate}
\end{proposition}
\proof Both graded rings in (1) are quotients of the Rees ring of the ideal $Q_2$, say $B[Q_2W] (\subset B[W])$. Let $L$ and $L'$ denote the corresponding homogeneous ideals of definition 
in $B[Q_2W]$, so $L\subset L'$. Let $[L]_n$ and $[L']_n$ denote the homogeneous components of degree $n$. Then, for each integer $n\geq 0$:
$$ [L]_n= Q_2^{n+1}+Q_1 Q_2^{n} (\subset Q_2^n) \mbox{  \ and \ } [L']_n= Q_2^{n+1}+(Q_1 \cap Q_2^{n})(\subset Q_2^n).$$

By assumption $Q_1 \cap Q_2^{n}$ is included in the integral closure of the ideal $Q_1 Q_2^{n}$ (in the ring $B$), and both are included in $Q_2^{n}$. So $\alpha \in Q_1 \cap Q_2^{n}\subset [L']_n$ fulfills a relation
$(\alpha)^N+a_1(\alpha)^{N-1}+\cdots + a_N=0$, with 
$a_i\in (\langle x_1, \dots, x_{e_1}\rangle \langle x_{e_2}, \dots, x_{d}\rangle^{n})^i$ in $S$. In particular $a_i\in (Q_1 Q_2^{n})^i\subset Q_1Q_2^{ni} \subset [L]_{ni} \subset Q_2^{ni}$, in the ring $B$, for $i=1, \dots, N$.
It follows that $L'$ is included in the integral closure of $L$ in the Rees ring $B[Q_2W]$.
In particular both $L$ and $L'$ have the same radical ideal in $B[Q_2W]$, which proves the claim in (1)
(the image of $L$ and $L'$ in 
$gr_{Q_2}(B)$ also have the same integral closure).

The claim in (2) is a consequence of the proof of (1) and Remark \ref{rk26}, 2).
\begin{parrafo}\label{Par28} We now address the question in (\ref{eq1}). Let $Y\subset X \subset W$ be inclusions of schemes, where $Y$ and $W$ are regular. Fix a closed point $x\in Y$ and assume that $x= Y\cap H$, where $H$ is a regular subscheme of $W$, and $Y$ and $H$ cut transversally and have complementary dimension. Let $W\leftarrow W'$ denote the blow up at $Y$, then the strict transform of $H$, say $H'$ induces a proper morphism, say $H\leftarrow H'$, which can be identified with the quadratic transformation of $H$ at $x$. In addition, there is a {\em natural identification} of the exceptional locus of $H\leftarrow H'$, with the fiber of $W\leftarrow W'$ over the closed point $x\in Y$.

It is natural to ask what remains of this identification of fibers when we replace $W$ by $X$ and $H$ by $\overline{X}=X\cap H$; when $X\stackrel{\pi}{\longleftarrow}  X'$ is the blow up at $Y$, and $\overline{X}\stackrel{\overline{\pi}}{\longleftarrow}  \overline{X}'$ is the blow up at $x$.

We study now this question in connection to equimultiplicity when $X$ has some additional properties. Let $(A,M)$ be an excellent pure dimensional ring of dimension $d$, let $P$ be a prime so that $(A/P, M/P)$ is regular (of dimension $e=d-ht(P)$). Fix
$x_{1}, \dots , x_e$ in $M$ which induce a regular system of parameters in $A/P$, and set $(\overline{A},\overline{M}=(A/\langle x_{1}, \dots , x_e \rangle
, M\overline{A})$. 

Note that $M=\langle x_1, \dots , x_e\rangle +P$, and  the previous question, reformulated in algebraic terms, leads us to the surjection
$ A/M\otimes_A gr_A(P) \to gr_{\overline{A}}(\overline{M}).$
 This is an isomorphism when $(A,M)$ is regular, which provides the natural identification of fibers in the regular case, and, in general, the surjection indicates that there is a closed embedding $\overline{\pi}^{-1}(x) \subset 
\pi^{-1}(x)$. The following Theorem \ref{TH211} characterizes equimultiplicity by the condition $\overline{\pi}^{-1}(x)_{red} = 
\pi^{-1}(x)_{red}$. In addition it will give us an extra criterion, used in our forthcoming discussion, which involves the tangent cone of $X$ at $x$.
\end{parrafo}
\begin{parrafo} \label{Par38}
Let $(A,M)$ be a local ring. Given a prime $P$ there is a naturally defined homomorphism $A/M\otimes_A gr_A(P)\to gr_A(M)$.
Assume, as above, that $(A/P, M/P)$ is regular of dimension, say $e$, that $x_{1}, \dots , x_e$ are elements in $M$ that induce a regular system of parameters in $A/P$, and hence $M=\langle x_{1}, \dots , x_e\rangle +P$. There is a surjective homomorphism, say 
$(A/M\otimes_A gr_A(P))[X_{1}, \dots, X_e]\to gr_A(M)\to 0$
which is non-canonical, defined by mapping $X_i$ to the class of $x_i$ in $M/M^2$, or say $$(A/M\otimes_A gr_A(P))\otimes_{A/M} gr_{A/P}(M/P) \to gr_A(M)\to 0.$$
\end{parrafo}
\begin{theorem}\label{TH211}{\em (Hironaka-Schickhoff)}\cite[Theorem 5]{Lipman2} Let $(A,M)$ be an excellent pure dimensional local ring of dimension $d$, with an infinite residue field. Let $P$ be a prime such that $A/P$ is regular, and fix $\{x_1, \dots, x_e\}$ as above. Set 

\begin{equation}\label{lec2}\xymatrix@R=-0.15pc@C=7pc{
H &  G\\
|| &  ||\\
(A/M\otimes_A gr_A(P))[X_{1},.,X_e]\ar[r]_{\Phi} \ar[dddd]^{\pi} & gr_A(M)\ar[dddd]^{\pi_1} \to 0\\
\\
\\
\\
A/M\otimes_A gr_A(P)  \ar[r]_{\Theta} & gr_{\overline{A}}(\overline{M}) \to 0\\
|| & ||
\\
\overline{H} & \overline{G}
}\end{equation} where
\begin{itemize}
\item $(\overline{A},
\overline{M})=(A/\langle x_1, \dots, x_e\rangle , M\overline{A});$
\item $ A/M\otimes_A gr_A(P))\to gr_A(M)$ is the natural map, and $\Phi$ maps $X_i$ on $x_i$;
\item $\Theta: A/M\otimes_R gr_A(P) \to gr_{\overline{A}}(\overline{M})$ is the natural surjection, as $P\overline{A}=\overline{M}$;
\item $\pi: (A/M\otimes_A gr_A(P))\otimes_{A/M} gr_{A/P}(M/P) =A/M\otimes_A gr_A(P)[X_1, \dots, X_e] \to A/M\otimes_A gr_A(P)$ is defined by setting $\pi(X_i)=0$, $i=1, \dots, e$.
\end{itemize}
Then the following conditions are equivalent:
\begin{enumerate}
\item The kernel of the surjection $\Phi$ is nilpotent.
\item $e_A(M)=e_{A_P}(PA_P)$.
\item The dimension of $\overline{A}$ is $d-e$ and the kernel of the surjection $\Theta$ is a nilpotent ideal.

\end{enumerate}
\end{theorem}
\proof $(1) \Rightarrow (2)$. If (1) holds, both rings have the same dimension. The dimension of $gr_A(M)=d$ (the dimension of $A$), and we conclude from this that the dimension of $A/M\otimes_A gr_A(P)$ is $d-e$, namely that $l_{(A,M)}(P)=ht(P)$, and the claim follows from our first characterization of equimultiplicity in Theorem \ref{TH14}.

$(2) \Rightarrow (3)$.  If 2) holds, then $l_{(A,M)} (P)=ht(P)$ (Theorem \ref{TH14}), so the dimension of the graded
ring $ A/M\otimes_A gr_A(P)$ is $ht(P)$. We assume that $A/M$ is an infinite field so one can choose elements $y_1, \dots, y_s$ in $P$, for $s=ht(P)$, which span a reduction of $P$. Recall that $M=\langle x_1, \dots , x_e\rangle +P$, and note that $M$ is
the integral closure of the ideal spanned by $y_1, \dots, y_s, x_1, \dots, x_e$. Here 
$s=ht(P)$, $e=dim(A/P)$, and $s+e=d$. So $\{y_1, \dots, y_s, x_1, \dots, x_e\}$ is a system of parameters. Finally, one can argue as in Theorem \ref{TH14}, so that (3) follows from Proposition \ref{P27}, (i).

$(3) \Rightarrow (1).$ 
Note that (3) implies that $dim( \overline{H})= dim( \overline{G})=d-e$. Recall that the dimension of a local ring equals the dimension of the graded ring for the latter equality.

As $dim (H)= dim( \overline{H})+e$, and the previous remark says that $d=dim(G)=dim \overline{G}+e$, we conclude that $d=dim(H)=dim(G)$.

We claim now that $G$ is pure dimensional. Let $\Spec(A) \leftarrow X$ denote the blow up at the closed point. Since $A$ is an excellent pure dimensional excellent scheme, the same holds at  $X$. Therefore all irreducible components of the closed subscheme of $X$ defined by the invertible ideal $M\calo_X$, namely $Proj(G)$, have the same dimension. Therefore the ring $G$ is also pure dimensional.

We finally turn to the study of $\Phi: (A/M\otimes_A gr_A(P))[X_{1}, \dots, X_e]\to gr_A(M)\to 0$ under the two assumptions obtained from (3): that both rings have the same dimension, and that $gr_A(M)$ is pure dimensional. This already shows that the ideal $ker(\Phi)$ is contained only in minimal primes.
To prove that $ker(\Phi)$ is nilpotent, we must show that it is contained in all the minimal primes of $H$. 
Let $q^*_i$, $i=1, \dots, m$ denote the minimal primes in 
$\overline{H}=A/M\otimes_A gr_A(P)$, and let 
$$q_i=q^*_iH, \ i=1, \dots, m,$$
denote the extension via the inclusion $\overline{H}=A/M\otimes_A gr_A(P)\subset H$. Clearly $q_i$, $i=1, \dots, m$, are the minimal primes in 
$H$.  We now fix an index $i_0$ and show that $ker (\Phi) \subset  q_{i_0}$.
The inclusions
$$ker (\Phi) \subset \pi^{-1}(ker (\Theta))=(ker (\Theta)) H+\langle X_1, \dots, X_e\rangle\subset 
q_{i_0}+\langle X_1, \dots, X_e\rangle, $$ follow from $ker (\Theta) \subset  q^*_{i_0}$, the definition of $\pi$, and the commutativity of the diagram.

As $\dim G= \dim H$ we note that the primes which contain $\ker  \Phi$, and are minimal with this condition, are also minimal primes in $H$. Therefore
the prime $q_{i_0}+\langle X_1, \dots, X_e\rangle$ contains a prime of the form $q_j=q^*_jH$, for some index $j$, say
 $ker (\Phi)\subset q_j \subset q_{i_0}+\langle X_1, \dots, X_e\rangle.$ So
$$(q_{i_0}+\langle X_1, \dots, X_e\rangle)\cap \overline{H} \supset q_j \cap \overline{H}=q_j^*$$
for some index $j$. On the other hand, as $H= \overline{H}[X_1, \dots, X_e]$, and $q_{i_0}= q^*_{i_0}H$, one checks that
$$(q_{i_0}+\langle X_1, \dots, X_e\rangle)\cap \overline{H} = q^*_{i_0} .$$
This shows that $i_0=j$, $q^*_{i_0}=q^*_{j} $ , so $ker (\Phi) \subset q_{i_0}$ for 
$i_0= 1, \dots, m$.
\endproof

\end{section}

\begin{section}{On the stability of transversality under blow-ups. }
\begin{parrafo}\label{condiciones}
We only consider here noetherian rings that are excellent, that contain a field,
 and with the property that all saturated chains of prime ideals have the same length.
 
 This class is closed by blow-ups, and if $B$ is in this class, the localization at any prime ideal, say $B_p$, is formally equidimensional (a requirement in Theorem \ref{Rees}, crucial  in what follows). 

As indicated in \ref{CincluidoB} one can usually assume that $B$ contains and is finite over a regular ring $S$, a fact that will show up in our discussion.
Set $L=B\otimes_SK$ where $K$ is the quotient field of $S$. In what follows set 
$n=\dim_K L$ (the generic rank of $B$ over $S$). Then Theorem \ref{MultForm} ensures that, for any prime $P\subset B$, the multiplicity $e_{B_P}(PB_P)$ is at most $n$. If $n=1$ then $S=B$ is regular.

We will assume that $n>1$, and $F_n(B)$ will denote the set of primes of multiplicity $n$ in $B$.
If $d$ denotes the Krull dimension of $B$, then $d$ is also the dimension of the regular ring $S$. By assumption, if $P'$ is a minimal prime in $B$, then $B/P'$ is $d$-dimensional and hence $P'\cap S=\{0\}$. By general properties of the multiplicity,
in order to study the set $F_n(B)$, we may assume that $B$ is included in $L$ (we may replace $B$ by its image in $L$), namely that non-zero elements in $S$ are non zero divisors in $B$. 
In fact, there is of course a surjective morphism 
\begin{equation}\label{411} B\to B'\subset B\otimes_SK
\end{equation}
by taking $B'$ as the image, so there is a closed immersion $\Spec(B')\subset \Spec(B)$. Since $B$ is pure dimensional, the kernel of $B\to B'$, say $J$, is the intersection of the p-primary components of the zero ideal in $B$ corresponding to the minimal prime ideals of $B$. In particular, $J$ is supported in dimension at most $d-1$, and one checks from this that
\begin{equation}\label{412}F_n(B)=F_n(B').
\end{equation}

The finite ring extension defines a finite morphism $\delta: \Spec(B) \to \Spec(S)$.
We say that an irreducible regular center $Y\subset \Spec(B)$ is {\em transversal} to the morphism if the multiplicity of $B_{Q}$ is the generic rank, where $Q$ denotes the generic point of $Y$. We show now that the generic rank is naturally compatible with blow ups at transversal centers. In addition, we show in Theorem \ref{16} that in order to study the highest multiplicity of $B$ it is possible to assume that it is complete intersection. 
Let us first recall some properties that hold in this setting.
\end{parrafo}
\begin{proposition} \label{prop1}
1) Let $S\subset B$ be a finite extension as above of generic rank n, then: 
\begin{itemize}
\item[1i)] The highest possible multiplicity along primes in $B$ is 
at most $n$. 
\item[1ii)] $\delta: \Spec(B) \to \Spec(S)$ maps $F_{n}(B)$ homeomorphically into its image $\delta(F_{n}(B))$.
\end{itemize}
2) Let $Q\subset B$ be a prime in $F_{n}(B)$, and assume that $B/Q$ is regular. Let $\p=Q\cap S$:
\begin{itemize}
\item[2i)] The natural inclusion $S/\p \subset B/Q$ is trivial, namely $S/\p = B/Q$, in particular both are regular rings.

\item [2ii)] There are elements $\theta_1 , \dots, \theta_m$ in $Q$ so that $B=S[\theta_1 , \dots, \theta_m]$, and $Q=\langle \p B, \theta_1 , \dots, \theta_m \rangle$. 
\item [2iii)] $Q$ is the integral closure of $\p B$ in $B$, in particular the elements $\theta_1 , \dots, \theta_m$ can be chosen in the integral closure of $\p B$.

\end{itemize}
 \end{proposition} 
 
  \proof The statements
   in 1) and also 2i) have been proved in Theorem \ref{T18}. 
 
 
 2ii) As $B$ is finite over $S$, $B=S[\theta'_1 , \dots, \theta'_m]$ for suitable elements $\theta'_1 , \dots, \theta'_m$ in $B$. Finally, as 
 $B/Q=S/\p$, one can find elements $\lambda_1, \dots, \lambda_m\in S$, so that $\theta=\theta'_i-\lambda_i\in Q$. Clearly 
 $B=S[\theta_1 , \dots, \theta_m]$, and $Q=\langle \p B, \theta_1 , \dots, \theta_m \rangle$.
 
2iii) It suffices to prove the claim at the localization $B_{\M}$ for every prime $\M$ containing $Q$.  
We claim that if $Q\in F_{n}(B)$, also $\M\in F_{n}(B)$. Firstly note that the multiplicity of $B$ at $\M$ is at most $n=\dim_K(B\otimes_SK)$. A theorem of Nagata ensures that $e_B(\M) \geq e_B(Q)$. His proof requires conditions on the local ring $B_{\M}$ which are fulfilled as we assume $B_{\M}$ to be excellent and pure dimensional (\cite[p. 153]{Nagata}).
So $\M\in F_{n}(B)$ for any prime $\M$ containing $Q$.

Set $\p=Q\cap S$, and $\m=\M\cap S$. Note that 
$S_{\m}/\p S_{\m}$ is regular.

Under these conditions $S_{\m}$ is regular and $\p$ is a regular prime, so one readily checks that $l(\p)=h(\p)$ at $S_{\m}$, namely that all fibers of the blow up at $\p$, over points in $V(\p) \subset \Spec(S)$, have the same dimension (see also \cite[Corollary 4.6, p.135]{Lipman2}).

As $\M$ is the only prime in $B$ dominating $\m$, $B_{\M}$ is finite over $S_{\m}$, so the blow up of $B_{\M}$ at $\p B_{\M}$ is finite over the blow up of $S_{\m}$ at $\p$. In particular 
the fibers over the closed points have the same dimension and therefore
$l(\p B_{\M})=h(\p B_{\M})$.

Note here that 
$Q=\sqrt{\p}$ in $B_{\M}$ and the conditions in \ref{C16} hold for $S_{\p}\subset B_{Q}$, so $\p B_P $ is a reduction of $P B_P$. Finally 2iii) follows from Böger's Theorem \ref{TB}.
\endproof

\begin{parrafo}\label{loc}{\em (On localizations and quotients.)} It is clear that  the generic rank of $S\subset B$ is not affected when $S$ is replaced by a localization. Moreover, the results in the previous Proposition \ref{prop1} are clearly compatible with localizations at $S$.
Our next result shows that there is a compatibility with quotients, say of $B$ by an ideal $J$. We  assume here that $B/J$ is pure dimensional so as to remain in the class \ref{condiciones}.

\begin{corollary}\label{cor44} Let the notation and conditions be as in \ref{prop1}. Assume that the generic rank of $B/J$ over $S$ is an integer $m\geq 1$ (in particular, we assume that $B$ and $B/J$ have the same dimension), then $F_n(B)\subset \Spec(B/J)$, and furthermore, $F_n(B) \subset F_m(B/J)$.
\end{corollary}
\proof As $m=\dim_K (B/J)\geq 1$, the finite morphism $\Spec(B/J) \to \Spec(S)$ is surjective. We may enlarge $J$ and assume that $B/J \subset K\otimes_S(B/J)$. Choose $\p \in \delta(F_n(B))$. There is a unique prime $P$ in $B$ dominating $\p$ (Corollary \ref{C16}), and $P\in F_n(B)$. Since $\Spec(B/J) \to \Spec(S)$ is surjective, $P$ must contain $J$. So it induces a prime, say $\overline{P}$ in $B/J$. Finally localize $S$ at $\p$ and check that as the conditions in 2) of Corollary \ref{C16} hold for $P$ at $B_P$, they also hold at $(B/J)_{\overline{P}}$.
\end{parrafo} 

\begin{parrafo}\label{NNN}
Fix $S\subset B$ subject to the conditions in \ref{condiciones}.
Set $B=S[\theta_1, \dots , \theta_m ]$. We will assume that $\dim_K( B\otimes_SK)=n>1$, and that 
\begin{equation}\label{inclusionn}B \subset B\otimes_SK.
\end{equation}
Let $f_1(Z), \ldots , f_m(Z) \in K[Z]$ denote the minimal polynomials of $\theta_1, \dots, \theta_m$ over the field $K$. In \cite{multiplicity}, Lemma 5.2, we show that if (\ref{inclusionn}) holds, then $f_1(Z), \ldots , f_m(Z) \in S[Z]$ (i.e., the minimal polynomials have coefficients in $S$). 

Let $d_i$ denote the degree of $f_i(Z)$, we may and will assume that $d_i\geq 2$, $i=1, \dots, m.$ Namely,
\begin{equation}\label{eq211}
B=S[\theta_1, \dots , \theta_m] \mbox{ and } d_i=\deg(f_i(Z))\geq2 \ , i=1, \dots ,m.
\end{equation}
The following result provides a description of the highest multiplicity locus.
\begin{theorem}\label{16} \cite[5.7]{multiplicity} Fix $B=S[\theta_1, \dots , \theta_m]$ 
as above. Then, if $ \dim_K(B\otimes_s K)=n $, a point ${\p}\in Spec(S)$ is the image of 
a point of multiplicity $n$ in $Spec(B)$, if and only if 
${\p}$ is the image of a point of multiplicity $d_i (\geq 2)$ in 
$Spec(S[Z]/\langle f_i(Z) \rangle)$, for every index $i=1,\dots , m$.
\end{theorem}
\end{parrafo}

\begin{parrafo}\label{centro}
Fix an inclusion $S\subset B$ of generic rank $n$ as in Prop \ref{prop1}, and a prime $Q$ in $F_n(B)$, so that $B/Q$ is regular. Consider now a prime $\M$ containing $Q$ (as in the proof 2iii) of Prop \ref{prop1}). Let $\m$, $\p$, in $S$, denote the primes dominated by $\M$ and $Q$ respectively. 
2iii) of Prop \ref{prop1} says that the natural inclusion of Rees rings
$S[\p W]\subset B[Q W]$ is a finite extension of graded rings. In particular the inclusion defines a finite morphism, say $X=\Proj (B[Q W])\to R=\Proj(S[q W])$.
Here $\Spec(S) \leftarrow R$ is the blow up of a regular scheme at a regular center. So given the previous setting, and after replacing $S$ by $S_f$ for a suitable $f\notin \m$, we may assume that there are elements $\{x_1, \dots, x_r, x_{r+1}, \dots, x_d \}$ in $S$,  which are a regular system of parameters at $S_\m$, and such that $\{x_1, \dots, x_r\}$ span the regular ideal $\p$. In this way one can cover $R$ by charts of the form $\Spec(S_t)$
\begin{equation}\label{cartas}
S_t=S\left[ \frac{x_1}{x_t}, \dots, \frac{x_r}{x_t}\right] (\subset K),  t=1, \dots, r.
\end{equation}

Fix the setting and notation as in (\ref{eq211}). Recall that we assume that $B/Q$ is regular, 
and the discussion in the proof of 2ii) shows that in this case one can choose elements, say 
$\lambda_1, \dots ,\lambda_m$ in $S$ (at least locally at any maximal ideal containing $Q$)
so that $\theta_i-\lambda_i \in Q$, for $i=1, \dots, m$. Note that the minimal polynomial 
of the element $\theta_i-\lambda_i $ can be extracted from the minimal polynomial of $\theta_i$ (from $f_i(Z)$) by a change of the variable $Z$ in $S[Z]$. Let
\begin{equation}\label{cambiov}
h_i(Z'_i),  \ i=1, \dots , m,
\end{equation}
denote the minimal polynomial of $\theta_i-\lambda_i $, so that 
$d_i= deg(f_i(Z))= deg(h_i(Z_i'))\geq 2$, and
\begin{equation}\label{eq212}
B=S[\theta_1-\lambda_1, \dots , \theta_m -\lambda_m].
\end{equation}

The following result is to be understood as a form of stability of the setting in Theorem \ref{16} (stability of transversality) when blowing up at smooth equimultiple centers.
\end{parrafo}

\begin{theorem}\label{Th120} \cite[Section 6.8]{multiplicity} Fix the notations and conditions as before, in particular assume that (\ref{inclusionn}) holds, and let $n$ denote the generic rank of $B$ over $S$ ($n=\mbox{dim}_K (B\otimes_SK)) $. Let $d=\dim B=\dim S$. Assume that $Y\subset \Spec(B)$ is regular irreducible with generic point $Q$, that it is included in the $F_n(B)\subset \Spec(B)$, and that $\dim Y<d$. Then a commutative diagram
\begin{equation}\label{edcdt} \xymatrix{
 Spec(B)\ar[d]^-{\delta} & & X\ar[d] ^-{\delta_1}\ar[ll]^-{\pi}\\
 Spec(S)&   & \ar[ll]^-{\pi'}R\\
 } 
\end{equation}
is given, where the horizontal morphisms are the blow ups at $Y$ and $\delta(Y)$ respectively, and $\delta_1$ is a finite dominant morphism. 
Moreover:

1)  Given a point $\m \in \delta(Y)$, and after taking a restriction of $\Spec(S)$ and $\delta: \Spec(B) \to \Spec(S)$ to a suitable affine neighborhood of $\m$, the blow-up $R$ can be covered by affine charts $U_t=\Spec(S_t)$, $t=1, \dots , r$, as in (\ref{cartas}), and $X$ by charts $V_t=\delta_1^{-1}(U_t)=\Spec(B_t)$, where
\begin{equation}\label{cartas1}
B_t=S_t\left[\frac{\theta_1-\lambda_1}{x_t}, \dots , \frac{\theta_M-\lambda_M}{x_t}\right].
\end{equation}

2) Non-zero elements in $S_t$ are non-zero divisors in $B_t$ and 
$K\otimes_{S_t} B_t=K\otimes_{S} B$ is the total quotient field of $B_t$. In particular $n=\mbox{dim}_K (B_t\otimes_{S_t}K) $, and the conditions in \ref{condiciones} hold.

3) The minimal polynomial of $\frac{\theta_i-\lambda_i}{x_t}$ is 
\begin{equation}\label{xyz11}
h'_i(W_i)=W_i^{d_i}+\frac{a_1^{(i)}}{x_t}W_i^{d_i-1}+ \cdots +\frac{a_{d_i}^{(i)}}{x_t^{d_i}}\in 
S_r[W_i],
\end{equation}
namely, the strict transform of that in (\ref{cambiov})(and $W_i$ denotes the strict transform of $Z_i'$).
\end{theorem}

\begin{parrafo}\label{laden} (On the proof of the Theorem). Let us outline the main ingredients for the proof which is very simple, and refer to \cite{multiplicity} for more details. To fix ideas take, as before, two primes 
$Q\subset \M$ in $B$, both of multiplicity $n$. Set $\p= Q\cap S$ and $\m=M\cap S$.
We assume here that $B/Q$ is a regular ring.  Proposition \ref{prop1} says that 
 that $\p$ is regular in $S_{\m}$  and $\p B$ is a reduction of $Q$. This latter condition 
 says that the natural inclusion
 $
 S_{\m}[\p W] \subset S_{\m}\otimes_SB[QW]
 $
is a finite extension. This finiteness already ensures the existence of the square diagrams in (\ref{edcdt}).
 The construction of the charts and the assertions 1), 2), and 3), are rather straightforward by checking the claims at each affine chart of the two blow ups (\ref{cartas1}).
\end{parrafo}

We now sketch a proof of the following theorem, also due to Dade (\cite{D}). Given $B$ as before, with maximum multiplicity $n$ at its prime ideals, our proof uses the existence of a finite extension $S\subset B$, where $S$ is a regular domain, and the generic rank of this extension is $n$.
This is not an additional condition for rings in the class which we are considering (see \ref{condiciones}), where one can argue locally. One produce such finite extension at the completion of the localization of $B$ at a prime.

\begin{theorem}\label{Dade} Let $X \stackrel{\pi_Y}{\longleftarrow} X'$ be the blow up at a regular center $Y$, equimultiple at $x\in Y$. Let $x'\in X'$ map to $x$, then 
$e(\calo_{X,x}) \geq e(\calo_{X',x'})$. In particular, if $Y$ is regular and equimultiple
$$\max Mult(X) \geq \max Mult(X').$$
\end{theorem}
\proof In our discussion we consider $X$ to be excellent and pure dimensional. Replace $X$ by $\Spec(B)$, an open neighborhood of $x$, and fix the notation as in the previous Theorem \ref{Th120} for a suitable inclusion $S\subset B$, of generic rank $n=e(\calo_{X,x})$. Note that for each index $t$, the generic rank of $B_t$ over $S_t$ is the same of that of $B$ over $S$. So the 
claim follows from 1i) in Proposition \ref{prop1}.
\begin{remark}\label{rk412}  Fix the setting as in Corollary \ref{cor44} for $\Spec(B/J) \subset \Spec(B)$. Fix $Y \subset F_n(B) (\subset F_m(B/J))$ as in Theorem \ref{Th120}, and note that the diagrams as in (\ref{edcdt}), for the blow ups of $\Spec(B)$ and $\Spec(B/J)$ respectively, induce 
an inclusion at the blow ups, say $X'_1 \subset X'$. Moreover, the previous conditions hold for this closed embedding, in particular $F_n(X') \subset F_m(X_1')$. 
\end{remark}
\end{section}
\begin{section}{Approaching a scheme by a complete intersection.}
We discuss here an interesting property of the multiplicity studied in \cite{multiplicity}. We show there that in order to reduce 
the multiplicity of a scheme by blow ups at equimultiple centers, it suffices to consider the case in which the scheme is a complete intersection. In fact, to a complete intersection with further additional properties; for example an equimultiple center will also be normally flat in the sense of Hironaka. 

Starting in \ref{rk415} we present some technical results which lead to our (inductive) Lemma \ref{417}, crucial for the proof of Main Theorem \ref{T3}. 
\begin{remark}\label{rk49} Set 
$B=S[\theta_1, \dots , \theta_m] \mbox{ and } d_i=\deg(f_i(Z))\geq2 ,  \ i=1, \dots ,m,$
as in (\ref{eq211}), where $S\subset B$ is of generic rank $n\geq 2$ (i.e., $\dim_K(B\otimes_SK)=n\geq 2$).

Let $S[V_1, \dots, V_m]$ be a polynomial ring over $S$ with $m$-variables, and consider now 
\begin{equation}\label{CI}B'=S[V_1, \dots, V_m]/\langle f_1(V_1), \dots, f_m(V_m) \rangle=
S[V_1]/\langle f_1(V_1) \rangle \otimes_S \dots \otimes_S S[V_m]/\langle f_m(V_m) \rangle.
\end{equation}
Note that $B'$ is a finite and free extension of $S$ of rank $D=d_1\cdot d_2 \cdot \dots \cdot d_m$. 

Set $B'=S[\overline{V}_1, \dots, \overline{V}_m]$ where $\overline{V}_i$ is the class of 
$V_i$, and note that $f_i(V_i)\in S[V_1, \dots, V_m]$ is the minimal polynomial of $\overline{V}_i$ over $S$. There is, in addition,  a natural surjection of finite $S$-algebras $$B'=S[\overline{V}_1, \dots, \overline{V}_m]\to B=S[\theta_1, \dots, \theta_m],$$ where each $\overline{V}_i$ maps 
to $\theta_i$, and both elements have the same minimal polynomial over $S$. Let 
$\delta: \Spec(B) \to \Spec(S)$ and $\delta ': \Spec(B') \to \Spec(S)$ be the corresponding finite morphisms, let $F_n(B)\subset \Spec(B) $ denote the points of multiplicity $n$, and 
$F_D(B')\subset \Spec(B') $ be those of multiplicity $D$. 
A consequence of Theorem \ref{16} is that 
\begin{equation}\label{ig1}
\delta(F_n(B))=\delta '(F_D(B')).
\end{equation}
Moreover, the natural inclusions $\Spec(B) \subset \Spec(B') \subset W=\Spec(S[V_1, \dots, V_m])$, where the latter is the regular scheme $W= \mathbb{A}_{S}^m$, together with the bijections given in Prop \ref{prop1} ensures that
\begin{equation}\label{ig2}
F_n(B)=F_D(B')  \subset W.
\end{equation}

Let $H_i$ be the hypersurface in $W$ defined by the polynomial $f_i(V_i) \in S[V_1, \dots, V_m]$, so
\begin{equation}\label{ig3}
\Spec(B')=H_1 \cap H_2 \cap \dots \cap H_m  \subset W,
\end{equation}
(scheme theoretical intersection).
Each hypersurface $H_i$ has at most multiplicy $d_i$, and 
\begin{equation}\label{ig4}
F_n(B)=F_D(B')=F_{d_1}(H_1) \cap \dots \cap F_{d_m}(H_m).
\end{equation}
Let $Y\subset \Spec(B)$ be a regular center included in $F_n(B)$ and consider the diagram
(\ref{edcdt})
defined by the blow ups.
Recall that the regular scheme $R$ is covered by charts given by rings
$S_t=S\left[\frac{x_1}{x_t}, \dots, \frac{x_r}{x_t}\right](\subset K),  t=1, \dots, r $ (see (\ref{cartas})),
and $\delta_1^{-1}\Spec(S_t)=\Spec(B_t)$, for $B_t=S_t\left[\frac{\theta_1-\lambda_1}{x_t}, \dots , \frac{\theta_m-\lambda_m}{x_t} \right]$ (\ref{cartas1}).
In addition, Theorem \ref{Th120} says that the minimal polynomial of $\frac{\theta_i-\lambda_i}{x_t}$ is the strict transform of $h_i$ in (\ref{cambiov}), which is obtained from $f_i$ by a change of variables (see (\ref{xyz11})). So the inclusion $\Spec(B)\subset H_1 \cap \dots \cap H_m$ 
induces the inclusion, say 
$$ X_1 \subset H'_1 \cap \dots \cap H'_m$$
after blowing up, where $H'_i$ is the strict transform of $H_i$ (see also Remark \ref{MN}). So again 
\begin{equation}\label{ig6}
F_n(X_1)=F_{d_1}(H'_1) \cap \dots \cap F_{d_m}(H'_m)  
\end{equation}
as the Theorem also provides a natural transformation of the data in (\ref{eq211}). 
%

The same argument, applied by iteration, shows that similar properties hold for any sequence of blow-ups at regular centers included in the set of points of multiplicity $n$. Namely a diagram, say
\begin{equation}
\label{gdigra}
\xymatrix@R=0pt@C=30pt{
X_0:=\Spec(B)\ar[dddd]^{\delta_0}     & \ar[l]
 X_1  \ar[dddd]^{\delta_1}   & \ldots  \ar[l] &   \ar[l] X_r \ar[dddd]^{\delta_r}  \\
\\
\\
\\
R_0:=\Spec(S)   & 
 \ar[l] R_1      & \ldots   \ar[l] &   \ar[l] R_r
}
\end{equation}
is defined, where the vertical maps  are finite morphisms that induce homeomorphisms, say 
$\delta_i: F_n(X_i)\to \delta_i(F_n(X_i))$ for each index $i=1, \dots ,r.$ In addition, for each $i$ there is an inclusion 
\begin{equation}\label{pol}X_i \subset X'_i=H^{(i)}_1 \cap \dots \cap H^{(i)}_m\subset W_i
\end{equation}
where each $H^{(i)}_j$ is a hypersurface in the regular scheme $W_i$, and $H^{(i)}_j$ is the strict transform of $H^{(i-1)}_j$. Moreover, for $i=0,1, \dots, r$:
\begin{equation}\label{ultima} F_n(X_i)=F_D(X'_i)=F_{d_1}(H^{(i)}_1) \cap \dots \cap F_{d_m}(H^{(i)}_m).
\end{equation}
\end{remark}
\begin{remark}\label{RA} Let $x\in X$ be a closed point of multiplicity $n$. If $X$ is a pure dimensional scheme of finite type over a perfect field, then at a suitable \'etale 
neighborhood of the point, say $\Spec(B)$, there is a regular subring $S$ so that $S\subset B$ is of generic rank 
$n$. When the characteristic is zero one can use the the previous results to produce a Rees algebra in $\Spec(S)$, so that if the lower row in (\ref{gdigra}) is a {\em resolution} of this algebra, then the upper row is a simplification of the singularity (see \cite{multiplicity}{ Section 7.1}, and apply maximal contact to eliminate the variables $\{X_1, \dots, X_m\}$ in(2)). In other words, in characteristic zero the simplification of the multiplicity locally at $x\in X$ is rephrased as a {\em resolution} of a Rees algebra
 on the scheme $\Spec(S)$, which is regular and of the same dimension as $X$. If we set as before 
 $B=S[\theta_1,\dots, \theta_m]$, then the Rees algebra in $\Spec(S)$ is constructed in terms of the coefficients of the minimal polynomial $f_i(V)$ (see \ref{NNN}). 
\color{black}
\end{remark}
\begin{remark}\label{MN} We claim that each $X'_i$ is the strict transform of $X'_{i-1}$ in (\ref{pol}). To prove this we show that $\Spec(B')$ is a complete intersection with very strong additional properties. In fact, (\ref{CI}) presents it as a fiber product of hypersurfaces which are finite over a fixed regular space $\Spec(S)$. 
Let $Y$ is an irreducible regular center included in $F_D(B')(\subset \Spec(S[V_1, \dots, V_m]))$, defined by a prime $P$ in $S[V_1, \dots, V_m]$. Set $\p=P \cap S$, which is regular in $S$. The previous discussion shows that for each index $i=1, \dots, m,$ and after a change of the variable $V_i$ in $S[V_i]$, we can assume that $V_i\in P$, and $f_i(V_i)\in P^{d_i}$. In other words, one can assume that $P=\langle \p, V_1, \dots, V_m\rangle$ and 
$f_i(V_i)\in \langle \p,V_i\rangle^{d_i}$. Identify, for simplicity, $V_i$ with $In_P(V_i)$, so that
$gr_{S[V_1, \dots, V_m]}(P)=(gr_S(\p))[V_1, \dots, V_m]$.

One can check now that each $In_P(f_i(V_i))$ is monic in $V_i$, of degree $d_i$, and that $\{ In_P(f_i(V_i)), i=1, \dots, m \}$ form a regular system on the graded ring $gr_{S[V_1, \dots, V_m]}(P)$.
This implies that 
\begin{equation} \label{eqmn}In_P(\langle f_1(V_1), \dots, f_m(V_m)\rangle) = \langle In_P(f_i(V_i)), \dots,In_P(f_m(V_m)) \rangle,
\end{equation}
in $gr_{S[V_1, \dots, V_m]}(P)$, so if $Q'$ denotes the class of $P$ in $B'$, 
\begin{equation}\label{eq444} gr_{B'}(Q')=gr_{S[V_1, \dots, V_m]}(P)/\langle In_P(f_1(V_1)), \dots, In_P(f_m(V_m))\rangle,
\end{equation}
and furthermore, 
\begin{equation}\label{eq444g} gr_{B'}(Q')=\otimes_i (gr_{S[V_i]}(P_i))/\langle In_{P_i}(f_i(V_i)),
\end{equation}
where $P_i= \langle \p, V_i\rangle$, $i=1, \dots, m$.

As $\Spec(B')$ is the closed subscheme in $\Spec(S[V_1, \dots, V_m])$ defined by the ideal spanned by $f_1(V_1), \dots , f_m(V_m)$, the equality (\ref{eqmn}) ensures that the blow up at $Y$ is defined 
by the strict transform of these $m$ hypersurfaces.  This proves our claim, namely that $X'_i$ is the strict transform of $X'_{i-1}$ in (\ref{pol}) for $i=1, \dots, r$.

 In addition, we claim that $Y (\subset F_D(B'))$ is a normally flat center at $\Spec(B')$ (i.e., that $gr_{B'}(Q')$ is flat over $B'/Q'=S/\p$). Equivalently, that the Hilbert-Samuel function of $\Spec(B')$ is constant along closed points in $Y$ (see also \ref{par52}).  To check this note that $gr_{S[V_i]}(P_i))/\langle In_{P_i}(f_i(V_i))$ is finite and free of rank $d_i$ over $S[V_i]/P_i=S/\p$, for $i=1, \dots, m$. So the claim follows from  (\ref{eq444g}).

Therefore, and as was indicated in  (\ref{forma10}) of the Introduction, (\ref{ultima}) can be refined by:
\begin{equation}\label{eqder}F_D(X'_i)=\Max \mult_{X'_i}=\Max HS_{X'_i}=F_{d_1}(H^{(i)}_1) \cap \dots \cap F_{d_m}(H^{(i)}_m).
\end{equation} \color{black}
\begin{parrafo}\label{54} In connection with the previous discussion let us mention 
that this is a particular case of Hironaka's notion of {\em idealistic presentation}. And one can show that he hypersurfaces defined by $f_1(V_1), \dots , f_m(V_m)$, together with the integers $d_i$, $i=1, \dots,m$, form an idealistic presentation of the complete intersection $\Spec(B')$ (but not of $\Spec(B)$!), as a subscheme of the regular scheme $\Spec(S[V_1, \dots, V_m])$,
in particular:
\begin{enumerate}
\item $(F_D(B')=) \Max \mult_{\Spec(B')}=\Max HS_{\Spec(B')}$;
\item $X'_i$ is the strict transform of $X'_{i-1}$ in (\ref{pol}) for $i=1, \dots, r$;
\item $(F_D(X_i')=) \Max \mult_{X'_i}=\Max HS_{X'_i}$.
\end{enumerate}

(see \cite[ (A), (B), (A'), (B') in p. 52, or Theorem 1, p. 100]{Hironaka77}. And this would also prove (\ref{eqder}).
%
%
%
\end{parrafo}
\end{remark}

\begin{parrafo}\label{EC} The rest of this section is devoted to Lemma \ref{417}, used in the proof of Theorem \ref{T3}. 

We now add an additional condition
to the data in (\ref{eq211}), which will be useful for the proof of Theorem \ref{T3}. Suppose now that $S=S_1[V]$, where $S_1$ is regular of dimension $d-1$, and fix a monic polynomial $g_b(V)=V^b+c_1V^{b-1}+ \cdots + c_b \in S_1[V]$. Let 
$H_b\subset \Spec(S)$ be the hypersurface defined by this polynomial. Note that 
$b$ is an upper bound for the multiplicity of $H_b$ at any point. Let $F_b(H_b)$ denote 
the (closed) set of points of multiplicity $b$, and assume that the following hold:
\begin{enumerate}
\item[(a)] $\delta(F_n(B)) \subset F_b(H_b)$ ($\subset \Spec(S)$);
\item[(b)] for {\em any} sequence as in (\ref{gdigra}), $\delta(F_n(X_i)) \subset F_b(H^{(i)}_b)$ ($\subset R_i$), for $0\leq i \leq r$, where $H_b^{(i)}$ denotes the strict transform of $H_b$.
\end{enumerate}
\end{parrafo}

\begin{remark}\label{rk415} Let us indicate that the additional condition in \ref{EC} occurs very naturally already in the study of the multiplicity of a hypersurface, and it will lead us to the construction of 
a new scheme in Lemma \ref{417}, which is a complete intersection
with properties similar to those of $B'$ in (\ref{CI}). 
Fix an index $i_0$, $1\leq i_0 \leq m$. 
Fix a point in $x\in W=\Spec(S[V_{i_0}])$ where the hypersurface defined by 
$f_{i_0}(V_{i_0})$ has multiplicity $d_{i_0}$. So $f_{i_0}(V_{i_0})$ has order $d_{i_0}$ at the local regular ring 
$\calo_{W,x}$, and the initial form $In_x(f_{i_0})$ is homogeneous of degree $d_{i_0}$ in the polynomial ring $gr_{\calo_{W,x}} (M_x)$. Assume that the residue field of $\calo_{W,x}$, say $k(x)$, is perfect. In this case Hironaka's $\tau$-invariant of this hypersurface at $x$ is the least number of variables needed to express $In_x(f_{i_0})$ in this graded ring. So $\tau(x) \geq 1$; and $\tau(x)=1$ if and only if $In_x(f_{i_0})$ is a power of a linear form.
On the other hand, if $\tau(x) \geq 2$ then at a suitable neighborhood one can set $S=S_1[V]$, where now $S_1$ is regular of dimension $d-1$ and $V$ is a variable over $S_1$, and find a monic polynomial $g_b(V)=V^b+c_1V^{b-1}+ \cdots + c_b \in (S=) S_1[V]$, so that $H_b\subset \Spec(S)$ fulfills the condition in \ref{EC}. A proof of this fact, where techniques of elimination are used, follows from \cite[Proposition 5.12]{hpositive}.

Let $\Spec(B') \subset \Spec(S[V_1, \dots, V_m])$ be given by (\ref{CI}), let $P$ be the prime in $S[V_1, \dots, V_m]$ corresponding to the point $x$, and $Q'=PB'$. Note that (\ref{eqmn}) and (\ref{eq444}) ensure that
\begin{equation}\label{eq444x} gr_{B_{Q'}'}(Q' B'_{Q'})=gr_{S[V_1, \dots, V_m]_P}(PS[V_1, \dots, V_m]_P)/\langle In_P(f_1(V_1)), \dots, In_P(f_m(V_m))\rangle.
\end{equation}
The previous discussion says that either the additional condition in \ref{EC} holds, or each initial form $In_P(f_i(V_i))$ is a power of a linear form. In this latter case 
the ring $(gr_{B_{Q'}'}(Q' B'_{Q'}))_{red}$ is regular.
\end{remark}

\begin{lemma}\label{417} Let the setting and notation be as in \ref{EC}. Set $\overline{B'}=B'/\langle g_b(V)\rangle,$ which is finite over $S_1$ (of dimension $d-1$). Let $s$ be the highest multiplicity at primes of $\overline{B'}$. Then

\begin{enumerate}
\item $F_n(B)=F_D(B')= F_s(\overline{B'})$ ($\subset \Spec(B')$).
\item For {\em any} sequence as in (\ref{gdigra}), $F_n(X_i)= F_D(X'_i)=F_s(\overline{X'}_i)$ ($\subset X'_i$), for $0\leq i \leq r$, where $\overline{X'}_i \subset X'_i$ denotes the strict transform of 
$\Spec(\overline{B'})$.
\item Any sequence $\overline{X'} \leftarrow \overline{X'} _1 \leftarrow \dots \leftarrow \overline{X'} _r$, defined by blowing up regular equimultiple centers included in the closed sets of points of multiplicity $s$, induces a sequence as in (2), where, in particular, $F_n(X_i)= F_s(\overline{X'}_i)$, $0\leq i \leq r$.
\end{enumerate}
\end{lemma}
\proof 
Consider $B'$ and $\overline{B'}$ as quotients of $S[V_1, \dots, V_m]$. $B'$ is defined by the equations 
$f_i(V_i)\in S[V_i]$, $i=1, \dots,m$, and these, together with $g_b(V)\in S=S_1[V]$ define
$\overline{B'}$. Choose a regular prime $P$ in $S[V_1, \dots, V_m]$ as in \ref{MN}, which induces a prime $Q'$ in $B'$, of multiplicity $D$, and $\overline{Q'}$ in  $\overline{B'}$.

Recall that we can assume that $V_i\in P$ and each $In_P(f_i(V_i))$ is monic in $V_i$, of degree $d_i$, and hence that $\{ In_P(f_1(V_1)), i=1, \dots, m \}$ form a regular system on the graded ring $gr_{S[V_1, \dots, V_m]}(P)$. The same argument applies for $g_b(V)\in S=S_1[V]$, so assume that $V\in P$ and that $In_P(g_b(V))$ is monic in $V$, in 
$gr_{S[V_1, \dots, V_m]}(P)=gr_{S_1[V, V_1, \dots, V_m]}(P)$.

We also conclude that $\{ In_P(f_1(V_1)), i=1, \dots, m \} \cup \{In_P(g_b(V))\}$ is a regular sequence.
This ensures that 
\begin{equation} \label{eqmncc}In_P(\langle f_1(V_1), \dots, f_m(V_m), g_b(V)\rangle) = \langle In_P(f_i(V_i)), \dots,In_P(f_m(V_m)), In_P(g_b(V)) \rangle \mbox{, so}
\end{equation}
\begin{equation}\label{eq444cc} gr_{\overline{B'}}(\overline{Q'})=gr_{S[V_1, \dots, V_m]}(P)/\langle In_P(f_1(V_1)), \dots, In_P(f_m(V_m)), In_P(g_b(V)) \rangle.
\end{equation}
Therefore the strict transform of $\Spec(\overline{B'})\subset \Spec(S[V_1, \dots, V_m])$
is defined by the strict transform of $\{ f_1(V_1), \dots, f_m(V_m), g_b(V)\}$. 
The Lemma follows now from the assumption in (\ref{EC}), and the commutative diagrams (\ref{gdigra}).

\end{section}

\begin{section}{Group schemes and the stratification of affine cones}

\begin{parrafo} The main result in this section is given in Theorem \ref{TH57} where we study the stratification defined by the multiplicity on an affine pure dimensional cone. This will lead us to the notion of group schemes,
also used in the formulation of Theorems \ref{T1} and \ref{T2}. We recall very briefly this concept and specify in \ref {propiedad} the properties used here. We refer to \cite[Chapter 2]{D-G} for this notion. 

Fix a commutative ring, $k$. In order to define a group (or monoid)  structure on a (affine) $k$-scheme $\Spec(B)$, it is convenient to view the scheme as a functor from the category of $k$-algebras to that of sets: $$G_B:\mathcal{A}_k \to {Ens},$$
namely, as the functor mapping the $k$-algebra $A$ to the set of homomorphisms of $k$-algebras $G_B(A)=Hom_{k-alg}(B,A)$. To this end we fix a finitely generated $k$-algebra $B$ together with a $k$-homomorphism, say $\Delta: B\to B\otimes_kB$. Then, for each $A\in \mathcal{A}_k$, and using the universal property of tensor products, we get a function
$$Hom_{k-alg}(B,A)\times Hom_{k-alg}(B,A)\to Hom_{k-alg}(B,A),$$
which is functorial in $A$. Roughly speaking, we say that the pair $(B, \Delta)$ defines an affine group (monoid) scheme when the previous function defines a group (monoid) structure on that set, and this structure is functorial. In such case
$G_B:\mathcal{A}_k \to \mathcal{G} \  (G_B:\mathcal{A}_k \to \mathcal{M})$
is a functor with values in the category of groups (monoids).

A first example arises when we fix a finitely generated module $M$ over a ring $k$, and consider the symmetric algebra 
$S_k[M]=B$. In this case, the universal properties of symmetric algebras says that we can also view $G_B(A)$ as the {\em group} of morphisms of $k$-modules, namely
$$G_B(A):=Hom_{k-alg}(B,A)=Hom_{k-mod}(M,A).$$
Note that the morphism of $k$-modules  $\delta: M \to M\oplus M$, $\delta(m)=(m,m)$ induces a homomorphism
$$\Delta: S_k[M] \to S_k[M\oplus M]=S_k[M]\otimes_kS_k[M],$$
of $k$-algebras. One checks that the pair $(S_k[M], \Delta)$ defines a group scheme, and each $A\in \mathcal{A}_k$ we recover the group structure of $Hom_{k-mod}(M,A).$

Suppose, finally, that we fix two $k$-algebras $B$, and $B'$, and two homomorphisms, say 
$\Delta: B \to B\otimes_kB$ and $\Theta: B' \to B\otimes_k B'$, and assume that the pair $(B,\Delta)$ defines a group or monoid structure on $G_B$. Note that $\Theta$ induces, for each $A\in  \mathcal{A}_k,$ a function
$$G_B(A)\times G_{B'}(A) \to G_{B'}(A).$$
which is functorial in $A$. We say that {\em the affine group scheme (monoid scheme) $G_{B}$ acts on $G_{B'}$} if, for each $A\in  \mathcal{A}_k,$ the previous map is an action of the group (monoid) $G_B(A)$ on the set 
 $G_{B'}(A)$.

We now discuss an example of particular interest in the study of affine cones (of graded $k$ algebras): Consider the affine monoid scheme $G_{k[X]}$, given by $(B=k[X], \Delta)$, where $\Delta: k[X] \to k[X]\otimes_kk[X]=k[X_1, X_2]$ maps $X$ to $X_1X_2$. Fix a $k$-algebra $B'$, and
a homomorphism $\Theta: B' \to B\otimes_k B'=B'[X]$. One can check that $\Theta$ defines an action 
of the affine monoid scheme $G_{k[X]}$ on $G_{B'}$ if and only if $\Theta: B' \to B'[X]$ defines a graded structure on $B'$ in a natural manner. By this we mean that a graded structure is defined 
on $B'$ so that $\Theta$ becomes a homomorphism of graded rings, where $B'[X]$ is endowed with the grading of the powers of $X$. In particular, if we fix a $k$-module $M$, the natural grading on the symmetric algebra $S_k[M]$ can be seen as a homomorphism $S_k[M] \to S_k[M]\otimes_kk[X]=S_k[M][X]$, which maps $m$ to $mX$ for each $m\in M$. One can finally check that this induces an action of the affine monoid scheme $G_{k[X]}$ on the affine scheme $G_{S_k[M]}$. Here $G_{k[X]}$ is to be understood as the scheme of scalars, acting on the abelian group scheme $G_{S_k[M]}$. 

 In our discussion we fix a field $k$, which will be perfect for most applications, and a finite dimensional $k$-vector space $M$. So $G_{S_k[M]}$ is a group scheme, with the natural action of the monoid $G_{k[X]}$ (with the natural grading on the symmetric algebra). This enables us to define the {\em vector space} $\mathbb{V}=\Spec(S_k[M])$ as this affine scheme considered together with these two operations. 
 
{\bf Invariant functions:} Let us mention that there is a natural notion of subspace, say $\mathbb{S}$, which is also a vector space.
In addition there is an action of this subgroup on $\mathbb{V}$: by $t_v: \mathbb{V} \to \mathbb{V}$, $t_v(w)=v+w$, for $v\in \mathbb{S}$.

There is also a naturally defined {\em subring of invariant functions}, say $I_{\mathbb{S}}\subset S_k[M],$ with the property that 
$\Spec(I_{\mathbb{S}})$, or say the functor $G_{I_{\mathbb{S}}}$, is identified with $\mathbb{V}/\mathbb{S}$ as group schemes and as a vector space.
In this construction, if $\mathbb{V}$ is a finite dimensional $k$-vector space defined by the polynomial ring $S_k[M]$, there is a natural identification of $M$ with the dual space $\mathbb{V}^{*}$.
In addition, to a subspace $\mathbb{S} \subset \mathbb{V}$ we assign the anulator $\mathbb{S}^{(0)} (\subset \mathbb{V}^{*}=M)$. The inclusion $\mathbb{S} \subset \mathbb{V}$ is given by the surjection 
$S_k[M] \to S_k[M]/\mathbb{S}^{(0)} S_k[M]$ . So $S_k[M]$ is a polynomial ring, and 
$\mathbb{S}^{(0)} S_k[M]$ is an ideal spanned by linear forms.

In this case the subring of functions of $\mathbb{V}$ which are invariant by the action of this subspace, acting by translations, is the symmetric algebra $S_k[\mathbb{S}^{(0)}]$, namely 
\begin{equation}\label{inclusion}I_{\mathbb{S}}=S_k[\mathbb{S}^{(0)}] \subset S_k[\mathbb{V}^{*}],
\end{equation}
and hence $I_{\mathbb{S}}$ is a graded $k$-algebra generated by linear equations, and it is also a polynomial ring. 

To be precise, functions in $I_{\mathbb{S}}$ are characterized as those functions on the scheme 
$\mathbb{V}$, invariant by all translations $t_v: \mathbb{V} \to \mathbb{V}$, $t_v(w)=v+w$, for $v\in \mathbb{S}$. Group schemes provide the right setting to make this concept precise, and to set $\Spec(I_{\mathbb{S}})=\mathbb{V}/\mathbb{S}$.

\end{parrafo}

\begin{parrafo}\label{propiedad}

An affine cone $C\subset \mathbb{V}$ is the closed subscheme defined by a homogeneous ideal $I(C)$ in $S_k[\mathbb{V}^{*}]$. We will say that the subspace 
$\mathbb{S}$ acts on $C$, when $t_v(C)=C$ for all translations with $v\in \mathbb{S}$, and
again, the language of group schemes provide a precise formulation. 
We will make use of the following properties for a cone $C\subset \mathbb{V}$ and a subspace $\mathbb{S}$ in $\mathbb{V}$:
\begin{enumerate}
\item There are (non-canonical) isomorphisms of schemes $(\mathbb{V}/\mathbb{S}) \times \mathbb{S} \to \mathbb{V}$.
\item $\mathbb{S}$ acts on $C$ if and only if there is a cone $C'\subset \mathbb{V}/\mathbb{S}=\Spec(I_{\mathbb{S}})$ so that (1) induces an isomorphism 
$C' \times \mathbb{S} \to C$.
\item If $\mathbb{S}_1$ and $\mathbb{S}_2$ are subspaces in $\mathbb{V}$ acting on $C$, then the subspace
$\mathbb{S}_1+\mathbb{S}_2$ acts on $C$.
\end{enumerate}

(3) follows from the definition of an action of a subspace (see \cite[ Def 5.2, p.29]{GiraudNotas}). Note that the inclusion in (\ref{inclusion}) is canonical, and the first two properties can be reformulated by:

\begin{enumerate}
\item[(1')] There are variables $V_1, \dots, V_r$ and (non-canonical) isomorphisms of graded $k$-algebras 
$$S_k[\mathbb{S}^{(0)}][V_1, \dots, V_r]\to S_k[\mathbb{V}^{*}]$$
extending the inclusion $S_k[\mathbb{S}^{(0)}] \subset S_k[\mathbb{V}^{*}]$.
\item[(2')]  $\mathbb{S}$ acts on $C$ if and only if $I(C)$ is the extension of an ideal 
in $I_{\mathbb{S}}=S_k[\mathbb{S}^{(0)}]$,  via (\ref{inclusion}). 
\end{enumerate}
1) follows from the previous discussion, setting $M=\mathbb{S}^{(0)}\oplus L$, where $L$ is a subspace with basis $V_1, \dots, V_r$. We refer here to \cite[Prop 5.4, p.32]{GiraudNotas} for the claims 2) or 2').
\end{parrafo}

\begin{parrafo}\label{51} Fix a polynomial ring $\P=k[X_1, \dots, X_r]$, over a field $k$, and let 
$\mathcal{A}$ be the quotient defined by a homogeneous ideal (i.e., let $C=\mathcal{A}$ be a graded algebra generated over $k$ by finitely many homogeneous elements of degree one). So
$\P=S_k[L]$ is the symmetric algebra over $k$ of the vector space of linear forms, say $L$, and
$\Spec(\mathcal{A})$ is a cone included in the vector space $\Spec(\P)$.

Let $\mathcal{A}_+$ denote the irrelevant ideal. This is a maximal ideal, rational over $k$, corresponding to the origin of the cone, say $\mathbb{O} \in C$.  Let $n$ denote the multiplicity of the cone at the origin, namely the multiplicity of $\mathcal{A}$ at the prime $\mathcal{A}_+$, and we want to study the set $n$ fold points $F_n(\mathcal{A})$.

The graded ring of $\mathcal{A}$ at the maximal ideal $\mathcal{A}_+$ can be identified with $\mathcal{A}$, namely 
\begin{equation}\label{A+} \mathcal{A}=gr_{\mathcal{A}}(\mathcal{A}_+).
\end{equation} We will always 
assume here that $\mathcal{A}$ is 
{\em pure dimensional}, say of dimension $d$. After taking a finite extension of $k$, in case it is not infinite, one can choose a polynomial subring, $k[Y_1, \dots, Y_d] \subset \mathcal{A}$, where $Y_1, \dots, Y_d$ are homogeneous of degree one, 
so that this is a finite extension of graded rings. There are many ways to choose such polynomial subrings, and the multiplicity of $\mathcal{A}$ at the prime $\mathcal{A}_+$ is given by the generic rank of this extension (see Corollary \ref{C16}).

As the ring is pure dimensional locally at 
$\mathcal{A}_+$, we may assume that $\mathcal{A}\subset \mathcal{A}\otimes_{k[Y_1, \dots, Y_d] }k(Y_1, \dots, Y_d)$.
In fact, if $\mathcal{A}'$ denotes the image of $\mathcal{A}$ in this total quotient ring, there is an identification, say 
$F_n(\mathcal{A})=F_n(\mathcal{A}')$ via $\Spec(\mathcal{A}') \subset \Spec(\mathcal{A})$ (see (\ref{412})).

%

\end{parrafo}

\begin{remark} \label{rk52} Choose $S=k[Y_1, \dots, Y_d] \subset \mathcal{A}$ as above, and assume that $\mathcal{A}\subset \mathcal{A}\otimes_{k[Y_1, \dots, Y_d] }k(Y_1, \dots, Y_d)$. Express $\mathcal{A}=S[\theta_1, \dots, \theta_m]$ choosing $\theta_i \notin S$, homogeneous of degree one.

Set $L=[\mathcal{A}]_1$, which we can assume to be free of rank $d+m$. Fix a basis 
\begin{equation}\label{grado1}\{Y_1, \dots, Y_d, V_1=\theta_1, \dots, V_m=\theta_m\}
\end{equation}
and set $\mathbb{V}=\Spec(S_k[L])$
where $S_k[L]$ is a polynomial ring over $k$. The surjection $S_k[L]\to \mathcal{A}$ induces an inclusion $\Spec(\mathcal{A}) \subset \mathbb{V}$. 
Let $F_i(Y_1, \dots, Y_d,V)$ denote the minimal polynomial of $\theta_i$ over $S$. As rings are graded $F_i$ is homogeneous in $k[Y_1, \dots, Y_d,V]$ of degree, say $d'_i$. Finally set 
$$ \mathcal{A}'=k[Y_1, \dots, Y_d, V_1, \dots, V_m]/ \langle F_1(Y_1, \dots, Y_d,V_1), \dots, F_M(Y_1, \dots, Y_d,V_m) \rangle.$$
The natural surjection, $ \mathcal{A}'\to  \mathcal{A}$ induces an isomorphism 
in degree one, say $[\mathcal{A}']_1= [\mathcal{A}]_1,$ and hence an inclusion 
$\Spec( \mathcal{A}) \subset \Spec( \mathcal{A}') \subset \mathbb{V}$.
According to the discussion is Section 4, and as a particular case of Theorem \ref{16}, the set of points in $\Spec(\mathcal{A})$ of multiplicity $n$, 
coincides with the set of points of multiplicity $D'=d'_1\cdot \dots \cdot d'_m$ in $\Spec(\mathcal{A}')$. Namely
\begin{equation}\label{for51}F_n(\mathcal{A}) =F_{D'}(\mathcal{A}') (\subset \mathbb{V}).
\end{equation}

\end{remark}

\begin{corollary}\label{lem51} Fix the setting and notation as above, and let $n$ denote the multiplicity of $\mathcal{A}$ at $\mathcal{A}_+$. If $k$ is perfect, then the set $F_n(\mathcal{A})$ 
of points in $\Spec(\mathcal{A})$ of multiplicity $n$, is a linear subspace.  
\end{corollary}
We refer to \cite[Ex 2.12 , III-25]{GiraudNotas}  for an example which shows that the statement does not hold over non-perfect fields.

\proof (of the corollary):
As $k$ is perfect, the points 
in $\Spec(k[Y_1, \dots, Y_d,V]/\langle F_i(Y_1, \dots, Y_d,V) \rangle)$ of multiplicity $d_i$ 
define a linear subspace, at least set theoretically. This can be checked by applying differential operators of degree $\leq d_i$ (see \cite[Lemme 1.2.7 , III-8]{GiraudNotas}). These give rise to an ideal spanned by lineal equation if the characteristic is zero, and to an ideal spanned by powers of linear equations if the field is perfect. In particular in both cases the set $F_{d_i}$, of points where the hypersurface 
has multiplicity $d_i$, is a linear subspace.
So if $k$ is perfect, $F_{D'}(\mathcal{A}')$ is an intersection of subspaces in $\mathbb{V}$.
\endproof
\begin{remark}\label{prop52}
Fix the notation and conditions as in \ref{51} (where $k$ is any field), then 
\begin{enumerate}
\item[i)] $\Proj(\mathcal{A})$ is pure dimensional, of dimension $d-1$;
\item[ii)] The multiplicity at points in $\Proj(\mathcal{A})$ is at most $n$;
\item[iii)] A point at $\Proj(\mathcal{A})$ has multiplicity $n$ if and only if it is defined by a homogeneous prime in $F_n(\mathcal{A})$.
\end{enumerate}

\end{remark}
\proof The assertion in i) follows from the assumptions on $\mathcal{A}$. The inclusion 
$k[Y_1, \dots, Y_d] \subset \mathcal{A}$ defines  $\pi: \Spec(\mathcal{A}) (\subset \mathbb{A}^m) \to \mathbb{A}^d=\Spec(k[Y_1, \dots, Y_d] ).$
Note firstly that the generic rank of $k[Y_1, \dots, Y_d] \subset \mathcal{A}$ is $n$ (\ref{51}). 
We also conclude from Proposition \ref{prop1} that $F_n(\mathcal{A})$
 maps homeomorphically to say $\pi(F_n)$ in 
$\mathbb{A}^d$. This one to one correspondence is compatible with localizations at
$k[Y_1, \dots, Y_d] $ (\ref{loc}), in particular for 
$k[Y_1, \dots, Y_d]_{Y_i} \subset \mathcal{A}_{Y_i}=\mathcal{A}\otimes_{k[Y_1, \dots, Y_d]}k[Y_1, \dots, Y_d]_{Y_i}.$

Set $S_i$ as the degree zero component of the $\mathbb{Z}$-graded ring $k[Y_1, \dots, Y_d]_{Y_i} $, and 
$A_i$ the degree zero component of  $\mathcal{A}_{Y_i} $. So
$k[Y_1, \dots, Y_d]_{Y_i} =S_i[Y_i, Y_i^{-1}]$ as graded rings, and $\mathcal{A}_{Y_i} =A_i[Y_i, Y_i^{-1}]$. One readily checks that $S_i\subset A_i$ is a finite extension 
of generic rank $n$. This proves (ii) and (iii).
\color{black}
\begin{theorem}\label{TH57} Let $C \subset \mathbb{V}$ be a pure dimensional cone over a perfect field.
 The highest multiplicity locus of $C$, say $\mathbb{S} (\subset \mathbb{V})$, is a subspace in $\mathbb{V}$, and it is the Hilbert Samuel stratum containing the origin of $C_{red}$.
\end{theorem}
\begin{corollary}\label{COR57} The highest multiplicity locus of a pure dimensional cone $C$ over a perfect field is a subspace, say $\mathbb{S}$ in  $\mathbb{V}$, 
and it is the biggest subspace of $\mathbb{V}$ acting on $C_{red}$ in the sense of group schemes. Equivalently, $\mathbb{S}$ is a subspace; it acts and it is the biggest subspace so that $C_{red}$ can be defined by functions of the quotient scheme $\mathbb{V}/\mathbb{S}$ (see \ref{propiedad}).
\end{corollary}
\begin{parrafo}\label{PTH57}(Proof of Theorem \ref{TH57}).
We have studied the highest multiplicity locus of a $k$-algebra $\mathcal{A}$. 
Here we view $\mathcal{A}$ together with it's graded structure. Let $L$ be, as before, the homogeneous component of degree one, so there is a surjection from the symmetric algebra, say $S_k[L]\to \mathcal{A}$, and hence, an embedding of cones 
\begin{equation}\Spec(\mathcal{A}) \subset \mathbb{V}=\Spec(S_k[L]).
\end{equation}
If $k$ is perfect, \ref{lem51} says that the highest multiplicity locus is a subspace, say $\mathbb{S}$, and $\mathbb{S} \subset \Spec(\mathcal{A}) \subset \mathbb{V}.$

Let $x\in \Spec(\mathcal{A})$ be the origin of the cone (corresponding to the maximal ideal $\mathcal{A}_+$).
Let $Q\subset \mathcal{A}$ denote the ideal defining this subspace, so $Q \subset \mathcal{A}_+$ is an inclusion of primes in $\mathcal{A}$. Applying Theorem \ref{TH211} to this setting, as $Y=\mathbb{S}$ is {\em regular and equimultiple} at the origin, it defines  
\begin{equation}\label{eq21}
(\mathcal{A}/\mathcal{A}_+\otimes_{\mathcal{A}} gr_{Q}(\mathcal{A}))
\otimes_{\mathcal{A}/\mathcal{A}_+}\mathcal{A}/\mathcal{A}_+[X_1, \dots,X_r] \to gr_{\mathcal{A}_+}(\mathcal{A})\to 0,
\end{equation} 
which is a surjection with a nilpotent kernel. Here $k(x):=\mathcal{A}/\mathcal{A}_+=k$, and $\mathcal{A}/Q$ is the polynomial ring $k[X_1, \dots, X_r]$. 
Let $K$ be the kernel of $S_k[L] \to \mathcal{A}$, and 
let $Q'\subset S_k[L]$ be the prime corresponding to $Q\subset \mathcal{A}$. So $\mathcal{A}/Q=S_k[L]/Q'=k[X_1, \dots, X_r]$. Consider the diagram 
\begin{equation}\label{edcdt45} \xymatrix{0\ar[d]&&0\ar[d]\\
K'\ar[d] \ar[rr]&&K=In_{S_k[L]_+}(K)\ar[d]\\
0\to (k(x)\otimes_{S_{k[L]}}gr_{S_{k[L]}}(Q'))\otimes_{k(x)} gr_{S_k[L]/Q'}((S_k[L]/Q')_+) \ar[d] \ar[rr]^{}&  &S_k[L]=gr_{S_{k[L]}}(S_k[L]_+) \to 0\ar[d] \\
 (k(x)\otimes_{\mathcal{A}}gr_{\mathcal{A}}(Q)) \otimes_{k(x)} gr_{S_k[L]/Q'}((S_k[L]/Q')_+) \ar[d] \ar[rr]^{}  &   & \mathcal{A}=gr_{\mathcal{A}}(\mathcal{A}_+)\ar[d]\to 0\\ 
 0&&0
 } 
\end{equation}
where equalities in the second column are as in (\ref{A+}), and the {\em isomorphism} in the middle is constructed from the natural morphism $gr_{S_k[L]}(Q') \to gr_{S_k[L]}(S_k[L]_+)$. Note that (\ref{A+}) also says that $$k(x)\otimes_{k(x)} gr_{S_k[L]/Q'}((S_k[L]/Q')_+) =S_k[L]/Q'=k[X_1, \dots, X_r],$$ 
and the lower row is (\ref{eq21}). 

Finally the equimultiplicity of $\mathcal{A}$
along $Q$ ensures that $(K'\subset) K \subset \sqrt{K'}$ (see \ref{TH211}). 

In this setting $K'$ is the extension of, say $K''=ker(k(x)\otimes_{S_{k[L]}}gr_{S_{k[L]}}(Q') \rightarrow 
k(x)\otimes_{\mathcal{A}}gr_{\mathcal{A}}(Q)))$, in particular $K'$ is spanned by elements in the subring $(k(x)\otimes_{S_{k[L]}}gr_{S_{k[L]}}(Q'))$.

We finally obtain, by general properties of group schemes:
\begin{enumerate}
\item The subspace $\mathbb{S}= \Spec(S_k[L]/Q')$ acts on 
$\mathbb{V}= \Spec(S_k[L])$ and the subring 
of invariant functions is $(k(x)\otimes_{S_{k[L]}}gr_{S_{k[L]}}(Q')) \subset S_k[L]$ (see (\ref{inclusion}) and \ref{propiedad}).
\item One can express the vector space as a direct sum $L=L_1\oplus L_2$ so that the previous inclusion is given by $(k(x)\otimes_{S_{k[L]}}gr_{S_{k[L]}}(Q'))=S_k[L_1] \subset S_k[L].$
\end{enumerate}
It follows from (1) and (2) that the ideal $\sqrt{K'}$ (in the ring of invariant functions $S_k[L_1]$) extends to a radical ideal in $S_k[L]$. In particular $\sqrt{K}$ can also be generated by invariant functions, so $\mathbb{S}= \Spec(S_k[L]/Q')$ acts on the affine cone $C_{red}=\Spec(\mathcal{A}_{red})$, and hence 
$C_{red}=C' \times \mathbb{S} $ for some 
$C'\subset \mathbb{V}/\mathbb{S}=\Spec(I_{\mathbb{S}})$ (\ref{propiedad}). As we assume that $\mathbb{S}$ is the stratum of the multiplicity through the origin of $C$, it is also a stratum of multiplicity of $C_{red}$, so $C_{red}=C' \times \mathbb{S} $, and we finally conclude from this that $\mathbb{S}$ is also the Hilbert-Samuel stratum of $C_{red}$ containing the origin.
\endproof

\end{parrafo}

\begin{section}{Proofs of Theorems \ref{T3}, \ref{T1} and \ref{T2}}
\end{section}

\begin{parrafo}\label{par52}{\em Factorization of cones and fibers of blow ups.}

We study the blow up of a scheme along a regular {\em equimultiple} center, and this leads to the study of cones, and factorization of cones in a vector spaces, which was partially  treated in the previous lines.
A vector space that arises naturally is the {\em Zariski tangent space} at $x\in X$, say $$T_{X,x}=\Spec(S_k[m/m^2]),$$ where $S_k[m/m^2]$ is the symmetric algebra,  $m$ is the maximal ideal at the point $x$, and $k=k(x)$ is the residue field. It is endowed with a structure of vector space as was mentioned in the previous section. 
The tangent cone of $X$ at $x$ is defined by $C_{X,x}=\Spec(gr_m(\calo_{X,x}))$, and the surjection of graded rings $S_k[m/m^2]\to gr_m(\calo_{X,x})$ induces the inclusion 
$C_{X,x}\subset T_{X,x}$.

We now discuss about ingredients of the blow up encoded in the tangent cone.
The projective scheme defined by $C_{X,x}$ is the fiber
over $x$ of the blow up at this point. In addition to this, we will see that $C_{X,x}$ also encodes information
of the blow up at suitable centers $Y (\subset X)$ which are regular at $x$. 
Given a regular subscheme $Y$ of $X$, defined by an ideal, say $Q$, a cone 
$C_{X,Y}= \Spec( gr_X(Q))$ is defined, where $gr_X(Q)=\oplus \ Q^n/Q^{n+1}$.
Let $X=\Spec(B)$, and let $x\in Y\subset X$ be defined by prime ideals 
$Q\subset M $ in $B$. Set $X\stackrel{\pi_Y}{\longleftarrow} X'$
the blow up at $Y$. Then 
$Proj(gr_B(Q))\subset X'$ is the closed subscheme defined by the invertible ideal $Q\calo_{X'}$, and the fiber at $x$ of this blow up is
\begin{equation}
\pi_Y^{-1}(x)=\Proj (k(x) \otimes_Bgr_B(Q)) (\subset X').
\end{equation}

Recall that the natural morphism $gr_B(Q) \to gr_M(B)$ induces 
$k(x) \otimes_Bgr_B(Q) \to gr_M(B),$ where $k(x)=B_M/MB_M$. We define the {\em conormal bundle} as the affine scheme
$$C_{X,Y,x}=\Spec(k(x) \otimes_Bgr_B(Q)),$$
and we study now conditions that enable us to relate $C_{X,Y,x}$ to 
$C_{X,x}$ (see (\ref{edcdt1}) and (\ref{535})).

It is convenient and handy to reformulate the discussion in \ref{PTH57}, particularly the diagram (\ref{edcdt45}) at the completion $\hat{B}_M$, since this ring can be viewed as a quotient of a ring of formal power series, say $(R,N)= k(x)[[y_1, \dots, y_e]]$, where $e=dim_{k(x)}(M/M^2)$. This is a local regular ring, and there is a natural identification $M/M^2=N/N^2$.
Note that $gr_R(N)$ is the symmetric algebra $S_k[M/M^2]$, or say 
$$ \Spec(gr_R(N))=C_{W,x}=T_{W,x}=T_{X,x}$$
where $W=\Spec(R)$, and hence $C_{W,x}$ is a vector space in the sense of group schemes.

A center $Y\subset X$, which is regular at $x\in Y$, induces, at the completion, a regular center, say $Y' \subset W$ and a regular prime $Q'$ in $(R,N)= k(x)[[y_1, \dots, y_e]]$.
After a change of coordinates we may assume that $Q'=\langle y_1, \dots, y_s \rangle$ for some $s\leq e$.
The homomorphism $gr_R(Q') \to gr_R(N)$ induces 
\begin{equation}\label{eqx} 0\to k(x) \otimes_Rgr_R(Q') \to gr_R(N),
\end{equation} which is an inclusion as the ring and the ideal are regular. Set, as before,
$C_{W,Y,x}=\Spec(k(x) \otimes_Rgr_R(Q'))$. Here $gr_R(N)=k(x)[Y_1, \dots, Y_e]$, where 
$Y_i$ is the initial form of $y_i$, and the image of $k(x) \otimes_Rgr_R(Q')$ is the subring
$k(x)[Y_1, \dots, Y_s]$.
So $T_{Y,x}$ is a subgroup (and furthermore, a subspace) of the vector space $T_{W,x}=T_{X,x}$, acting by translations, and the ring of invariant functions on $T_{W,x}$ 
$k(x)[Y_1, \dots, Y_s]$. Group schemes enable us to identify (canonically):
$ C_{W,Y,x}=T_{W,x}/T_{Y,x}  \mbox{ (see \ref{propiedad})}.$
Recall that $C_{X,x}\subset T_{V,x}=T_{W,x}$, and the subspace $T_{Y,x}$ {\em acts on the cone} $C_{X,x}$ when this cone can be defined by equations in the subring of invariant function. Namely, when the ideal $I(C_{X,x})$ can be spanned by equations in the image of the injective homomorphism (\ref{eqx}).

Let us repeat (essentially) the diagram (\ref{edcdt45}) in this context, where $Q$ is regular and equimultiple in $B$.
Set $(R', N')=(R/Q', N/Q')$. So $gr_{R'}(N')$ is a polynomial ring, and (\ref{eqx}) can be extended in many ways to an {\em isomorphism} $\Phi'$, defining  a diagram of graded rings and ideals:

\begin{equation}\label{edcdt1}\xymatrix@R=-0.15pc@C=7pc{
0\ar[ddd]  & 0\ar[ddd] \\
& \\
& \\
& \\
K\ar[ddd] \ar[r] &  In(J)\ar[ddd] \\
&\\
&\\
 &  \\
0\to (k(x)\otimes_Rgr_R(Q'))\otimes_{k(x)} gr_{R'}(N') \ar[r]_{\Phi'} \ar[dddd]^{\pi} & gr_R(N) \ar[dddd]^{\pi_1} \to 0\\
\\
\\
\\
(k(x)\otimes_Bgr_B(Q)) \otimes_{k(x)} gr_{R'}(N')\ar[ddd]   \ar[r]_{\Phi} & gr_{B}(M)\ar[ddd]  \to 0\\
& \\
& \\
& \\
\\
0 & 0
}\end{equation}
and recall that $K$ is the ideal spanned by, say $K'=ker( (k(x)\otimes_Rgr_R(Q'))  \to (k(x)\otimes_Bgr_B(Q))).$
 We finally reformulate the discussion \ref{PTH57} in terms of factorization of cones by saying that this square diagram of graded rings induce:
 \begin{enumerate}
 \item[A)] An {\em isomorphism} at the homogeneous component of degree one, and hence a a non canonical diagram of schemes included in the Zariski tangent cone $ T_{X,x} $, say
 \begin{equation}\label{edcdt1} \xymatrix{
C_{W,Y,x}\times T_{Y,x}  &= &  T_{X,x} \\
 C_{X,Y,x}\times T_{Y,x} \ar[u]  &   & C_{X,x}\ar[ll]^{\Phi^*} \ar[u]\\
 } 
\end{equation}
where the vertical morphisms are given, essentially, by the {\em closed immersions} $C_{X,Y,x} \subset C_{W,Y,x}$
and $C_{X,x}\subset T_{X,x}$ (recall that $\Spec(gr_R(N))=T_{W,x}=T_{X,x}$). 

\item[B)] Canonical identifications  $C_{W,Y,x}=T_{X,x}/T_{Y,x}$, $C_{X,Y,x}=C_{X,x}/T_{Y,x}$ if $Y$ is normally flat at $x$, and $(C_{X,Y,x})_{red}=(C_{X,x})_{red}/T_{Y,x}
$ when $Y$ is equimultiple at $x$.
\end{enumerate}

In fact, in (\cite{Hironaka64},Theorem 2, p. 195) Hironaka proves that normal flatness along $Y$ occurs when $\Phi^*$ is a scheme theoretical identification
(when $\Phi$ is an isomorphism), whereas equimultiplicity 
along $Y$ 
occurs when 
$$0 \to ((k(x)\otimes_Bgr_B(Q)))_{red} \otimes_{k(x)} gr_{R'}(N')\to (gr_B(M))_{red}\to 0,$$
is an isomorphism.
%
%
or say, if and only if $\Phi^*$ induces an identification $(C_{X,x})_{red} = (C_{X,Y,x})_{red}\times T_{Y,x},$
as subschemes in the tangent space. We now repeat (for convenience) and prove  Theorem \ref{T1}.

\begin{theorem}
\label{PT1} 
Fix a point $x\in X$ and let $C_{X,x}\subset T_{X,x}$ be the inclusion of the tangent cone in the Zariski tangent space. Assume that the residue field $k(x)$ is perfect, and let $\mathbb{S}$ denote the stratum of (highest) multiplicity through the origin of the cone 
$C_{X,x}$.
\begin{enumerate}
\item $\mathbb{S}$ is a subspace in $T_{X,x}$. It acts, and it is the biggest subspace acting on $(C_{X,x})_{red}$. 

\item Let $Y\subset X$ be regular and equimultiple at $x$, then the subspace 
$T_{Y,x}$ is included in $\mathbb{S}$, and hence it also acts on the subscheme $(C_{X,x})_{red}$, and  $(C_{X,x})_{red}/T_{Y,x}=(C_{X,Y,x})_{red}.$
\item Let $X\stackrel{\pi_Y}{\longleftarrow} X'$ be the blow up at $Y$, as above. For any $x'\in \pi^{-1}(x) \subset X'$, 
$e(\calo_{X', x'}) \leq e(\calo_{X, x})$, and if the equality holds then
\begin{equation}\label{721}x' \in \Proj (\mathbb{S}/T_{Y,x})\subset \pi^{-1}(x)=\Proj(C_{X,Y,x}),
\end{equation}
where the inclusion is obtained from the action of $T_{Y,x}$ on $\mathbb{S}\subset (C_{X,x})_{red}$ in $T_{X,x}$.
\end{enumerate}
\end{theorem}

\proof
(1) and (2) where discussed in Theorem \ref{TH57}, Corollary \ref{COR57}, and in \ref{par52}.

(3) The inequality $e(\calo_{X', x'}) \leq e(\calo_{X, x})$ is a result of Dade, presented here in \ref{Dade}. So we prove now
(\ref{721}). Let 
$Q\subset \calo_{X,x}$ be the prime defining the regular center $Y$. In what follows set $B=\hat{\calo}_{X,x}$, the completion. Note that $k(x)\otimes_{\calo_{X,x}} gr_{\calo_{X,x}}(Q)=k(x)\otimes_{B} gr_{B}(QB)$. This enables us to assume that $X=\Spec(B)$, and that $X\stackrel{\pi_Y}{\longleftarrow} X'$ is the blow up at $QB$. Here $\mathbb{S}$ is also the stratum of highest multiplicity in the reduced scheme $(C_{X,x})_{red}$ (see (\ref{MULTRED})), and the discussion at the end of \ref{PTH57} (Proof of Theorem \ref{TH57}),
shows that 
\begin{equation}\label{??}\mathbb{S}/T_{Y,x}
\end{equation}
is the stratum of highest multiplicity in $(C_{X,Y, x})_{red}$ (or equivalently, of $C_{X,Y, x}$).

Set $Q'=QB$, and let $B[Q' W]$ be the corresponding Rees ring so $X'= Proj(B[Q' W])$. We must prove that if $x' \in X'$ in (3) is such that $e(\calo_{X', x'}) =e(\calo_{X, x})$, $x'$ corresponds to a homogeneous prime in $k(x)\otimes_{B} gr_{B}(QB)$ which is in the set of points of highest multiplicity of this graded ring (namely in $\mathbb{S}/T_{Y,x})$). We shall prove this in the next lines, and hence we prove (3). The argument will rely on the equivalence of (2) and (3) of Theorem \ref{TH211}. Before doing so, we discuss some preliminary technical aspects.

%
%
%

After replacing $x\in X$ by a suitable \'etale neighborhood in case the residue field is finite
(see \ref{finito}), one can assume that the conditions in Lemma \ref{L13}
hold at the completion $(B,M)$: Namely, although the residue field might not be infinite, one can assume that there is a system of parameters $\{y_1, \dots, y_d\}$, that span a reduction of $M$, and that the extension 
$S=k(x)[[y_1, \dots, y_d]]\subset B$ is finite with generic rank $[L,K]=e(\calo_{X,x}),$ where $K$ and $L$ are the corresponding quotient fields; and that $Q'$, is the integral closure of $\langle y_1, \dots, y_s \rangle B$ for some $s\leq d$. Note that the class of $\{y_{s+1}, \dots, y_d \}$ is a regular system of parameters at $B/Q'$, as they span a reduction of the maximal ideal, and the maximal ideal of a regular local ring does not admit a proper reduction. 

Set $(\overline{B}, \overline{M})=(B/\langle y_{s+1}, \dots, y_d \rangle , M\overline{B})$, which is finite over $S/\langle y_{s+1}, \dots, y_d \rangle=\overline{S}=k(x)[[y_1, \dots, y_s]]$. The previous construction ensures that the 
class of $y_1, \dots, y_s$, span a reduction of $\overline{M}$.

Set $\overline{K}=k(x)((y_1, \dots, y_s))$ (the quotient field of $\overline{S}$). After replacing $\overline{B}$ by its image in 
$\overline{L}=\overline{K}\otimes_{k(x)[[y_1, \dots, y_s]]} \overline{B}$ we may assume that 
the multiplicity of $(\overline{B}, \overline{M})$ is the generic rank $[\overline{L}, \overline{K}]$.
Let $N$ and $\overline{N}$ be the maximal ideals in $S$ and $\overline{S}$ respectively, and consider the diagram of Rees rings 
\begin{equation}\label{edtf} \xymatrix{
 S\oplus N \oplus N^2\oplus \cdots \ar[rr]^-{} \ar[d] ^-{}& & B\oplus Q' \oplus Q'^2\oplus \cdots\ar[d] ^-{}\\
 \overline{S} \oplus \overline{N} \oplus \overline{N}^2\oplus \cdots \ar[rr]^-{}&   & 
 \overline{B} \oplus \overline{M} \oplus \overline{M}^2\oplus \cdots\\
 } 
\end{equation}
where the horizontal morphisms are finite extensions and the vertical ones are surjective. It induces
\begin{equation}\label{edtf} \xymatrix{
 Z_1=\Proj (S\oplus N \oplus N^2\oplus )& & \ar[ll]^-{}X'=\Proj(B\oplus Q' \oplus Q'^2\oplus \cdots) \\
  \ar[u] ^-{}\overline{Z_1} =\Proj(\overline{S} \oplus \overline{N} \oplus \overline{N}^2\oplus \cdots)&   & 
 \ar[u] ^-{}\ar[ll]^-{}\overline{X}'=\Proj( \overline{B} \oplus \overline{M} \oplus \overline{M}^2\oplus \cdots)\\
 } 
\end{equation}
where the horizontal morphisms are finite, and the vertical ones are closed immersions.
Set $\Spec( \overline{B})= \overline{X}\subset X=\Spec(B)$ and $\Spec( \overline{S})= \overline{Z}\subset Z=\Spec(S)$. So 
\begin{enumerate}
\item[(a)] 
 $X=\Spec(B)$ is finite over $Z=\Spec(S)$, and $\overline{X}$ is finite over $\overline{Z}$.
\item[(b)] $\overline{X}'\subset X'$ and both are finite over the regular schemes $\overline{Z}_1\subset Z_1$ respectively.
\end{enumerate}
Fix $x'\in X' $ as in (3), so that:
\begin{enumerate}
\item[(i)] $x'$ maps to the closed point of $\Spec(B)$; and
\item[(ii)] $e(\calo_{X', x'})=e_B(M)$. 
\end{enumerate}

We claim that if $e(\calo_{X', x'}) =e(\calo_{X, x})$, then $x'$ corresponds to a homogeneous prime in $k(x)\otimes_{B} gr_{B}(QB)$ of highest multiplicity in this ring. The equivalence of (2) and (3) in Theorem \ref{TH211}
says that the natural surjection $k(x)\otimes_{B} gr_{B}(QB)\to gr_{\overline{B}}(\overline{M})$ induces an isomorphism of the {\em reduced} graded rings, so the spectrum of both graded rings share the same stratum of points of highest multiplicity (see (\ref{MULTRED})), namely
$\mathbb{S}/T_{Y,x}$ (\ref{??}). 

Therefore (i) says that $x'\in \overline{X}' (\subset X')$, and we claim that (ii) implies that:

\begin{enumerate}
\item[(ii')] $e(\calo_{\overline{X}', x'})=e_{\overline{B}}(\overline{M})$
$(=[\overline{L}:\overline{K}])$,
\end{enumerate}
where $\overline{X}\leftarrow \overline{X}'$ is now the blow up of 
$\Spec(\overline{B})$ at the origin. So (i) and (ii') would imply that $x'$ corresponds with a homogeneous prime in $\Spec(gr_{\overline{B}}(\overline{M}))$ which is in the stratum of the multiplicity containing 
the origin. Equivalently, that this homogeneous prime is in $\mathbb{S}/T_{Y,x} \subset \Spec(gr_{\overline{B}}(\overline{M}))$ (see Remark \ref{kr34}), and this would prove our claim in (\ref{721}).

We finally proof that (i) and (ii) implies (ii'). Let $z_1\in Z_1$ denote the image of $x'\in {X}'$. Note that the total quotient ring of 
$\calo_{X', x'}$ is $L$ (the total quotient ring of $B$), and the total quotient ring of 
$\calo_{Z_1, z_1}$ is $K$ (the quotient field of $S$). So the generic rank of the finite 
extension $\calo_{Z_1}\subset \calo_{X_1}$ is $e(\calo_{X, x})=e(\calo_{X', x'})$. This shows that $x'$ must be the only point mapping to $z_1$, so $\calo_{Z_1, z_1}\subset \calo_{X', x'}$ is finite. In fact Proposition \ref{prop1} shows that the maximal ideal of the local regular ring $\calo_{Z_1, z_1}$ spans a reduction of the maximal ideal of $\calo_{X', x'}$, and both have the same residue field.
As  $\calo_{\overline{X}', x'}$ is a quotient of $\calo_{X_1, x_1}$ and $\calo_{\overline{Z}_1, z_1}$ is a quotient of $\calo_{Z_1, z_1}$, one checks that 
$\calo_{\overline{Z}_1, z_1} \subset \calo_{\overline{X}', x'}$ inherits the previous properties. Namely, that this latter inclusion is finite, the maximal ideal of the local regular ring $\calo_{\overline{Z}_1, z_1}$ spans a reduction of the maximal ideal of $\calo_{\overline{X}', x'}$, and both have the same residue field.
The equivalence in Corollary \ref{C16} ensures finally that the multiplicity of $\calo_{\overline{X}', x'}$ is given by the generic rank over $\calo_{\overline{Z}_1, z_1}$, namely by $[\overline{L}:\overline{K}]$, and hence that $e(\calo_{\overline{X}', x'})=e_{\overline{B}}(\overline{M})$.
\endproof

\end{parrafo}

\begin{parrafo}\label{PT2}(Proof of Theorem \ref{T2}.)
Set $S\subset B$, $S\subset B'$, $B'\to B\to 0$, and $\x \in \Spec(B)( \subset \Spec(B'))$ with multiplicity 
$n$ and $D$ respectively, as in in Remark \ref{rk49}. Let $\M'$, $\M$, and $\mathcal N$
be the prime ideals in $B'$, $B$, and $S$, respectively, corresponding to $x$. So $\Spec(B')$ is a complete intersection, and as was indicated in \ref{MN}, the claims in (\ref{forma0}) and (\ref{forma1}) of our introduction follow from (\ref{ig6}), (\ref{ultima}).

Recall that all the schemes considered here are pure dimensional. Let $d$ be the dimension of the local ring $B_{\M}$ (of $B'_{\M'}$). There is a surjection $gr_{B'}(\M')\to gr_{B}(\M)$, inducing an inclusion $C_{X,x} \subset C_{X',x}$. 

 Let $\mathbb{S}_{HS}' \subset \mathbb{S}_{mult}'(\subset C_{X',x}=\Spec(gr_{B'}(\M'))$ be the inclusion of the stratum through the origin of the Hilbert Samuel function and of the multiplicity respectively. 
 Let $\mathbb{S}(\subset C_{X,x}=\Spec(gr_{B}(\M))$ be the stratum of the multiplicity.
 
We first note that there are inclusions $(\mathbb{S}_{HS}'\subset )\mathbb{S}_{mult}' \subset \mathbb{S}$. This follows by applying Corollary \ref{cor44} to this setting, as both graded rings
in the surjection 
$gr_{B'}(\M')\to gr_{B}(\M)$ are finite over $gr_S(\mathcal N)$, and hence of the same dimension.
%

In particular $\dim(\mathbb{S}_{HS}') \leq \dim(\mathbb{S}),$
and finally Theorem \ref{T2} is a corollary of (\ref{eqder}) 
and some well known properties of the $\tau$ invariant of a hypersurface mentioned in Remark \ref{rk415}. More generally, our Theorem \ref{T2} follows now from a result of Hironaka:
\end{parrafo}

\begin{theorem}\label{T22} Assume that the residue field $k(x)$ is perfect at a point $x\in X'$ and let $\mathbb{S}'$ denote the stratum of the Hilbert Samuel function  through the origin of the cone 
$C_{X',x}$. Assume, for simplicity that 
$x\in \Max HS_{X'} $, and let $D'(x)$ be the dimension of the subspace $\mathbb{S}'\subset T_{X',x}= T_{X,x}.$
\begin{enumerate}
\item $D'(x)$ is an upper bound 
of the local dimension of the closed set $\Max HS_{X'} $ at $x$.
\item For any sequence $X'\leftarrow X'_1 \leftarrow \dots \leftarrow X'_r$, of blow ups at normally flat centers, and given closed points $x_i\in X'_i$, $0\leq i \leq r$ so that $x_{i+1}$ maps to $x_i$, and 
$x_0=x$, then $$HS(\calo_{X'_i, x_i})\geq HS(\calo_{X'_{i+1}, x_{i+1}})$$ for 
$i=0, \dots , r-1$, and if $HS(\calo_{X'_r, x_r})=HS(\calo_{X', x})$, then $D'(x)$ is an upper bound 
of the local dimension of the closed set $\Max HS_{X'_r} $ at $x_r$.
\end{enumerate}
\end{theorem}
 This result concerns the Hilbert Samuel stratum. It follows from \cite[Part 4 of Chapter III.]{Hironaka64}, where 
the behavior of the $\tau$-invariant is studied in the formal setting. 
It is also treated in \cite{Hironaka77} in the setting of \'etale topology, where the Theorem is a consequence of 
the notion of Local Distinguished Presentations in Definition 6, p. 82, and of their Transformation in p. 94. 
\begin{parrafo}\label{PT3}(Proof of Theorem \ref{T3}). 
Fix again the setting and notation as in \ref{PT2}, where all schemes considered are pure dimensional. Let $d$ be the dimension of the local ring $B_{\M}$ (of $B'_{\M'}$). So the surjection $gr_{B'_{\M'}}(\M')\to gr_{B_{\M}}(\M)$, induces an inclusion $C_{X,x} \subset C_{X',x}$ of affine schemes of pure dimension $d$. An affine cone 
 is regular if and only if it is regular at the origin, in which case it is also 
irreducible. In particular, if $(C_{X',x})_{red}$ is regular, both cones 
$C_{X,x}$ and $C_{X',x}$ have the same underlying topological space, or say 
$(C_{X,x})_{red}=(C_{X',x})_{red}$, and hence $(C_{X,x})_{red}$ is also regular.
More generally, as $d$ is the dimension at $x\in X'$, and hence of 
$C_{X',x}$, the stratum of highest multiplicity in the cone, namely 
$\mathbb{S}' (\subset C_{X',x}$), is of dimension at most $d$. Here $B'$ be an affine chart
and the setting  of \ref{PT2} provides an inclusion $\Spec(B') \subset \Spec(S[V_1, \dots, V_m])$ as in \ref{rk49}. After replacing the regular ring $S$ by $S_{\mathcal N}$, for $\mathcal N= S\cap \M '$, $S_{\mathcal N}$ has dimension $d$. The discussion in Remark \ref{rk415}, and the presentation of $gr_{B'_{\M '}}(\M ' B'_{\M '})$ in (\ref{eq444x}), show 
that the inclusion $C_{X',x} \subset 
T_{X,x}$ is given by a surjection $gr_{S_{\mathcal{N}}[V_1, \dots, V_m]_P}(PS[V_1, \dots, V_m]_P)\to gr_{B'_{\M '}}(\M ' B'_{\M '}))$, where $P=\langle \mathcal{N}, V_1, \dots, V_m\rangle$.
Consider 
$$T=gr_{S_{\mathcal{N}}[V_1, \dots, V_m]_P}(PS[V_1, \dots, V_m]_P)
=gr_{S_{\mathcal{N}}}(\mathcal{N}S_{\mathcal{N}})[V_1, \dots, V_m],$$ (where we identify 
$V_i$ with it's initial form).
In Remark \ref{rk415} we also show that each $In_P(f_i(V_i))$ involves only $V_i$ (i.e., $In_P(f_i(V_i))\in gr_{S_{\mathcal{N}}}(\mathcal{N}S_{\mathcal{N}})[V_i]\subset T$). In addition, $gr_{S_{\mathcal{N}}}(\mathcal{N}S_{\mathcal{N}})$ is a polynomial ring in $d$ varables, and 
$T$ is a polynomial ring in $d+m$ variables over the perfect field 
$k(x)$. Finally, Remark \ref{rk415}  also stated that the least number of variables required to define the 
 equations $\{ In_P(f_1(V_1), \dots , In_P(f_m(V_m)\}$ in
 the polynomial ring $T$ is exactly $m$ (i.e., the dimension of $\mathbb{S}'$ is $d$), only 
 when each 
$In_P(f_i(V_i))$ is a power of a linear form in $T$. 
The proof of the Main Theorem \ref{T3} follows now from Lemma \ref{417}.

\end{parrafo}

\end{section}


\begin{thebibliography}{www}


\bibitem{Abad} C. Abad, 
{ `On the Highest Multiplicity Locus of Algebraic Varieties and Rees algebras'}.
Journal of Algebra.

\bibitem{Aab} S.S. Abhyankar,
{\em Ramification theoretic methods in Algebraic Geometry}, Annals
of Mathematics Studies  No. {43}, Princeton University Press (1959).
Princeton New Jersey.

%



%
%
%
%
\bibitem{BeV1} A. Benito, O. E. Villamayor U., `Monoidal transforms and invariants of singularities in positive characteristic'. 
{\em Compositio  Mathematica} Volume 149, no 8, (2013) 1267-1311.
\bibitem{BeV2} ------ `Techniques for the study of singularities with applications to resolution of 2-dimensional schemes´. {\em Math. Ann.} Volume 353, Number 3,(2012) 1037-1068 
\bibitem{BB} B.M. Bennett, `On the characteristic functions of a local ring´. {\em Ann. of Math}, 91,{(1970)}, 25-87.

%
%
\bibitem{BM} E. Bierstone, P. Milman,  `Canonical
 desingularization
in characteristic zero by blowing-up the maxima
strata of a local
invariant', {\em Inv. Math.},  128 (1997) 207-302.
%
%


\bibitem{Indiana} 
A. Bravo, M.L Garc\'{\i}a-Escamilla, O.E. Villamayor U., {`On Rees algebras and 
invariants for singularities over perfect fields'}. {\em Indiana University Mathematics Journal}. 
Issue 3 of Volume 61, (2012) 1201-1251.

\bibitem{BrV} A. Bravo, O. E. Villamayor U., `Singularities in positive characteristic, stratification and simplification of the singular locus', {\em Adv. in Math.}, 224 (4) (2010) 1349-1418. 

\bibitem{BrV1} ------ Notes on Resolution of Singularities: `On the behavior of the multiplicity on schemes:
stratification and blow ups'. 81-207. The Resolution of Singular Algebraic Varieties, Clay Institute Math. Proceedings, vol 20, D.Ellwood, H.Hauser,S.Mori, and J. Schicho Editors - AMS-CMI, 2014.
%
\bibitem{CJS} V. Cossart, U. Jannsen, S. Saito, `Canonical embedded and non-embedded resolution of singularities for excellent two-dimensional schemes'. Prepint: arXiv:0905.2191. 

\bibitem{CP1} V. Cossart, O. Piltant, `Resolution of singularities of threefolds in positive characteristic. I. Reduction to local uniformization on Artin-Schreier and purely inseparable coverings',  {\em J. Algebra}, 320 (3) (2008)   1051-1082. 

\bibitem{CP2} ------ `Resolution of singularities of threefolds in positive charac. II.',  {\em J. Algebra}, 321 (7) (2009) 1836-1976.

\bibitem{Cut} S. D. Cutkosky, `Resolution of singularities of 3-folds in positive characteristic.',  {\em Amer. J. Math.}, 131 (2009), no. 1, 59-127. 
%

\bibitem{Cut2} ------  {\em Resolution of Singularities} (Graduate Studies in Mathematics, Vol. 63).

\bibitem{D} E. C. Dade, {\em Multiplicity and monoidal transformations.} Thesis, Princeton University, 1960.

\bibitem{Dietel} B. Dietel. {\em A refinement of Hironaka?s additive group schemes for an extended invariant.} (Thesis) Regensburg 2014

\bibitem{D-G} M. Demazure, P. Gabriel, {\em Introduction to algebraic geometry and algebraic groups. } Mathematics Studies 39, North-Holland, 1980. 
%

\bibitem{EncVil97:Tirol}
 S.~Encinas and O. E. Villamayor U.,
 \newblock `A Course on Constructive Desingularization and Equivariance',
 \newblock  {\em Resolution of Singularities. A research textbook  in tribute to Oscar Zariski} (eds H. Hauser, J. Lipman, F. Oort, A. Quir\'os),
 Progr. in Math.  181 (Birkh\"auser, Basel, 2000) pp. 147-227.

\bibitem{GK} T. Gaffney, S. Kleiman, Specialization of integral dependence for modules,
Invent. Math.137 (1999) 541?574.

%
%
%
%

\bibitem{Giraud1975} J.~Giraud,
\newblock `Contact maximal en
caract\'eristique positive',
 \newblock {\em Ann. S. de l'Ec. N.}, 8
 (1975) 201-234.
 
 \bibitem{GiraudNotas} ------
\newblock `Etude locale des singularit\'es ',
 \newblock {\em Cours de 3\'eme Cycle, Univ. Paris XI,} U.E.R. Math. 91-Orsay, 26 (1971-72), (mimeographed).



\bibitem{EGAIV}{A. Grothendieck and J. Dieudonn\'e}, 
	\textit {\'El\'ements de G\'eom\'etrie Alg\'ebrique},
	\textit  {IV \'Etude locale des sch\'emas et des morphismes de sch\'emas}, {Publications Math\'ematiques} I.H.E.S. No. 24 (1965) and No. 28 (1966).
	 
	

\bibitem{Ha} H, Hauser, `The Hironaka Theorem on resolution of singularities'. Bull. Amer. Math. Soc. 40 (2003), 323-403.
\bibitem{Ha1} ------  `On the problem of resolution of singularities in positive characteristic (Or: A proof
we are still waiting for)'. Bull. Amer. Math. Soc. 47 (2010), 130.

\bibitem{Herrmann} M. Herrmann, S. Ikeda, U. Orbanz, {\em Equimultiplicity and Blowing up}, Springer-Verlag, 
Berlin, 1982, 200.232.

\bibitem{HerrmannAL} M. Herrmann et al.,  `Theorie der normalen Flachheit' , Leipzig. Teubner 1977.

\bibitem{HerrmannOr} M. Herrmann, U. Orbanz, {\em Between equimultiplicity and normal flatness}, Algebraic Geometry, Proceedings La Rabida 1981. Lecture Notes in Mathematics 961. Springer-Verlag, 
Berlin, 1988.


\bibitem{Hironaka64}  H.~Hironaka, `Resolution of singularities
 of
an algebraic variety over a field of characteristic zero I-II',
 {\em Ann.
Math.},  79 (1964) 109-326.
 \bibitem{HironakaS}------
`Normal cones in analytic Whitney stratifications'.
 Publ. IHES 36 (1969), 127-138.

 
 \bibitem{Hironaka777} ------
`Introduction to the theory of infinitely near singular points'.
 Memorias de Matematica del Instituto Jorge Juan. Vol. 28.
 Concejo Superior de Investigaciones Cient\'ificas. Madrid, 1974. 

 


%
\bibitem{Hironaka70} ------ 'Additive groups associated with points of a projective space',
\emph{Ann. of Math.} 92 (2) (1970), 327-334.
%
%
%
%
%
%

%
%


\bibitem{Hironaka77} ------
`Idealistic exponent of a
 singularity', {\em Algebraic Geometry},  The
John Hopkins centennial
 lectures, Baltimore, John Hopkins University Press
(1977), 52-125.





\bibitem{kaw} H. Kawanoue, `Toward resolution of singularities
over a field of positive characteristic I'.{\em Publ. Res. Inst. Math. Sci.}, 
43 (2007) 819-909. 

\bibitem{kaM} H. Kawanoue,  K. Matsuki, `Toward resolution of singularities over a field of positive characteristic' (the idealistic 
filtration program) Part II. Basic invariants associated to the idealistic filtration and their properties. Publ. Res. Inst. Math. Sci., 46(2):359{422}, 2010.




\bibitem{Lipman} J. Lipman, `Stable ideals and Arf rings', {\em Amer. J.  Math.}  93 (1971) 649-685. 



\bibitem{Lipman1} ------ `Relative Lipschitz-saturation', {\em Amer. J.  Math.}  97 (1975) 791-813. 

\bibitem{Lipman2} ------ `Equimultiplicity, reduction, and blowing up', {\em Commutative Algebra.} Editor: R.N.Draper.
Lecture Notes in Pure and Appl. Math. 68
Marcel Dekker, New York (1982) 111-147.

\bibitem{Mat} H.~Matsumura.
 {\em Commutative algebra},  {\em
Mathematics
    Lecture Note Series, 56},
 Benjamin/Cummings Publishing
Company, Inc., 2nd ed. edition, 1980.

\bibitem{Mat1} ------
 {\em Commutative ring theory},  {\em
Cambridge studies in advanced math., 8},
 Cambridge Univ. Press, 1989.

\bibitem{Nagata} M. Nagata, {\em Local Rings}, John Wiley and Sons, New York, 1962.

\bibitem{Nagata1}------`The theory of multiplicity in general local rings'  {\em Proceedings of the International Symposium. Tokyo-Nikko 1955. Scientific Council of Japan, Tokio. } (1956) 191-226.

%
\bibitem{Nob} A. Nobile, `Equivalence and resolution of singularities'. Journal of Algebra, Vol. {420}, (2014) 161--185.
\bibitem{NR} D. G. Northcott and D. Rees, `Reduction of ideals in local rings'  {\em Proc. Camb. Philos. Soc.} 50 (1954) 145-158.


\bibitem{Oda1987} T. Oda, `Infinitely very near singular points'  {\em Complex analytic
singularities},  Adv. Studies in Pure Math. {\bf8}  (North-Holland, 1987) 363--404.

\bibitem{Oda2} ------- `Hironaka's additive group scheme', in Number Theory, AlgebraicGeometry and Commutative Algebra in honor of Y. Akizuki (Y. Kusunoki et al., eds.), Kinokuniya, Tokyo, 1973, 181-219.
\bibitem{Orbanz} U. Orbanz, `Multiplicities and Hilbert functions under blowing up,'  {\em Man. Math.} 36 (1981), 179-186.

\bibitem{Pomerol} M.J. Pomerol, `Sur la strate de Samuel du sommet d'un cône en caractéristique positive'. Bull. Sc. math., 2e série 98, 1974, p. 173-182.
\bibitem{R} M. Raynaud, \textit{Anneaux Locaux Hens\'eliens}, 
Lectures Notes in Mathematics 169. Springer-Verlag, 
Berlin, Heidelberg, New York, (1970).

\bibitem{Rees}  D. Rees, `Transforms of local rings and a theorem of multiplicity,'  {\em Proc. Camb. Philos. Soc.} 57 (1961) 8-17.
%
\bibitem{Samuel}  P.Samuel, `La notion de multiplicit\'e 
en Algebra et en G\'eom\'etrie alg\'ebrique'.  {\em Journal de Math., tome XXX}  (1951) 159-274.
\bibitem{Sch} W. Schickhoff, 
`Whitneysche Tangentenkegel, Multiplizitätsvarietäten und Äquisingularitätstheorie für Ramissche Räume'. Univ., Diss.--Osnabrück, 1976.


\bibitem{BS}  B. Singh, `Effect of a permissible blowing up on local Hilbert functions,'  {\em Inventiones Math.,} 26 (1974) 201-212.

\bibitem{Villa89} O. E. Villamayor U., `Constructiveness of Hironaka's resolution', {\em Ann. Scient. Ec. N.
Sup.} 22 (1989) 1-32.


\bibitem{Villa92}
------`Patching local uniformizations', {\em Ann.
Scient. Ec. Norm. Sup.},  25 (1992), 629-677.
%
%


\bibitem{hpositive}
------ {`Hypersurface singularities in positive
characteristic'}, {\em Adv.  in Math.} 213 (2007) 687-733.



%
\bibitem{multiplicity} ------{Equimultiplicity, algebraic elimination, and blowing-up'}.
{\em Adv. in Math.} 262 (2014), 313-369. 
%
\bibitem{ZS} O. Zariski, P. Samuel, {\em Commutative Algebra,}Vol.2 Van Nostrand, 
Princeton, 1960.


 \bibitem{WL} J. Wlodarczyk, Simple Hironaka resolution in characteristic zero, {\em  J. Amer. Math. Soc.}  {\bf 18}  no. 4  (2005), 779--822.


\end{thebibliography}
\end{document}